\documentclass[11pt,a4paper]{article}
\usepackage{amsmath,amssymb}
\newcommand{\R}{\mathbb{R}}
\newcommand{\Z}{\mathbb{Z}}

\def\cC{{\mathcal C}}

\def\cS{{\mathcal S}}

\newcommand{\e}{\varepsilon}
\renewcommand{\div}{{\rm div}\,}

\newcommand{\Sup}{\displaystyle \sup}

\newcommand{\Sum}{\displaystyle \sum}

\def\ddl{\dot \Delta_l}

\def\ddq{\dot \Delta_q}
\def\tilde{\widetilde}
\def\hat{\widehat}
\newcommand{\va}{\varphi}
\newcommand{\D}{\Delta}

\newcommand{\n}{\nabla}
\newcommand{\N}{\frac{N}{2}}

\newcommand{\fd}{\frac{N}{2}}
\newcommand{\fdp}{\frac{N}{p}}
\newcommand{\NN}{\frac{N}{p}}
\newcommand{\p}{\partial}

\newcommand{\de}{\delta}

\newcommand{\dq}{\delta q}

\newcommand{\ub}{\overline{u}}

\newtheorem{theorem}{Theorem}
\newtheorem{lem}{Lemma}

\newtheorem{prop}{Proposition}
\newtheorem{proposition}{Proposition}
\newtheorem{defi}{Definition}

\newtheorem{remarka}{Remark}

\title{New formulation of the compressible Navier-Stokes equations and parabolicity of the density}

\author{Boris Haspot  \thanks{Ceremade UMR CNRS 7534
Universit\'e de Paris  Dauphine,
Place du Mar\'echal DeLattre De Tassigny
75775 PARIS CEDEX 16 , haspot@ceremade.dauphine.fr }}
\date{}
\begin{document}
\maketitle

\begin{abstract}
In this paper we give a new formulation of the compressible Navier-Stokes by introducing an suitable effective velocity $v=u+\n\va(\rho)$ provided that the viscosity coefficients verify the algebraic relation of \cite{BD}. We give in particular a very simple proof of the entropy discovered in \cite{BD}, in addition our argument show why the algebraic relation of \cite{BD} appears naturally. More precisely the system reads in a very surprising way as two parabolic equation on the density $\rho$ and the vorticity ${\rm curl}v$, and as a transport equation on the divergence ${\rm div}v$. We show  the existence of strong solution with large initial data in finite time when $(\rho_0-1)\in B^{\NN}_{p,1}$. A remarkable feature of this solution is the regularizing effects on the density. We extend this result to the case of global strong solution with small initial data.
\end{abstract}
\section{Introduction}
The motion of a general barotropic compressible fluid is described by the following system:
\begin{equation}
\begin{cases}
\begin{aligned}
&\p_{t}\rho+{\rm div}(\rho u)=0,\\
&\p_{t}(\rho u)+{\rm div}(\rho u\otimes u)-{\rm div}(2\mu(\rho)D(u))-\n(\lambda(\rho){\rm div} u)+\n P(\rho)=0,\\
&(\rho,u)_{/t=0}=(\rho_{0},u_{0}).
\end{aligned}
\end{cases}
\label{0.1}
\end{equation}
Here $u=u(t,x)\in\R^{N}$ stands for the velocity field and $\rho=\rho(t,x)\in\R^{+}$ is the density. As usual, $D(u)$ stands for the strain tensor, with $D(u)=\frac{1}{2}\left(^t\nabla u+\nabla u\right)$.

The pressure $P$ is a suitable smooth function depending on the density $\rho$.
We denote by $\lambda$ and $\mu$ the two viscosity coefficients of the fluid,
which are also assumed to depend on the density and which verify some parabolic conditions for the momentum equation $\mu>0$ and $\lambda+2\mu>0$ (in the physical cases the viscosity coefficients verify $\lambda+\frac{2\mu}{N}>0$ which is a particular case of the previous assumption).
We supplement the problem with initial condition $(\rho_{0},u_{0})$.
Throughout the paper, we assume that the space variable $x\in\R^{d}$ or to the periodic
box $\mathbb{T}^{N}_{a}$ with period $a_{i}$, in the i-th direction. We restrict ourselves to the case $N\geq2$. In the sequel, we shall assume that $\mu$ and $\lambda$ are regular functions and verify the following algebraic relation introduced in \cite{BD} (see also \cite{MV}):
\begin{equation}
\lambda(\rho)=2\rho \mu'(\rho)-2\mu(\rho).
\label{1.2}
\end{equation}
We define $\va(\rho)$ by $\va'(\rho)=\frac{2\mu'(\rho)}{\rho}$ (in the sequel $\va(\rho)$ will be well defined since we shall assume no vacuum on the density). 
As in \cite{cras,CPAM1,CPAM} setting $v=u+\n\va(\rho)$
the system (\ref{0.1}) becomes (see the appendix for more details):
\begin{equation}
\begin{cases}
\begin{aligned}
&\p_{t}\rho-2\D\mu(\rho)+{\rm div}(\rho v)=0,\\
&\rho \p_{t}v+\rho u\cdot\n v-{\rm div}(\mu(\rho){\rm curl}v)+\n P(\rho)=0,\\
&(\rho,u)_{/t=0}=(\rho_{0},u_{0}).
\end{aligned}
\end{cases}
\label{0.1a}
\end{equation}
In passing we observe that this new formulation of the compressible Navier-Stokes equation allows us very easily to obtain the entropy discovered in \cite{BD}. Indeed multiplying the momentum equation of (\ref{0.1a}) by $v$ we have using the fact that ${\rm curl}v={\rm curl}u$:
\begin{equation}
\begin{aligned}
&\int_{\R^N}\big(\rho(t,x)|v|^2(t,x)+(\Pi(\rho)(t,x)-\Pi(\bar{\rho})) \big)dx+\frac{1}{2}\int^t_0\int_{\R^N}\mu(\rho)|{\rm curl}u|^2 (t,x)dx dt\\
&+\int^t_0\int_{\R^N}\n P(\rho)\cdot\n\va(\rho)(t,x)dx dt\leq \int_{\R^N}\big(\rho_0(x)|v_0|^2(x)+\Pi(\rho_0)(x) \big)dx,
\end{aligned}
\label{entrop}
\end{equation}
with $\bar{\rho}>0$, $\Pi(s)=s(\int^s_{\bar{\rho}}\frac{P(z)}{z^2}dz-\frac{P(\bar{\rho})}{\bar{\rho}})$ and $p(\rho)=\rho^\gamma$ with $\gamma>1$. (\ref{entrop}) is exactly the entropy discovered in \cite{BD} (see also \cite{MV}).\\
Assuming now that $\rho>0$ and dividing by $\rho$ the momentum equation we obtain the following system:
\begin{equation}
\begin{cases}
\begin{aligned}
&\p_{t}\rho-2\D\mu(\rho)+{\rm div}(\rho v)=0,\\
& \p_{t}v+ u\cdot\n v-\frac{\mu(\rho)}{\rho}{\rm div}{\rm curl}v-\frac{1}{2}\n \va(\rho)\cdot{\rm curl}v+\n F(\rho)=0,
\end{aligned}
\end{cases}
\label{0.1ab}
\end{equation}
with $F'(\rho)=\frac{P'(\rho)}{\rho}$.  When we apply the ${\rm curl}$ and the ${\rm div}$ to the momentum equation, we obtain since
$\D{\rm curl} v={\rm curl}\,{\rm div}\,{\rm curl}v$ and ${\rm div}\,{\rm div}\,{\rm curl}v=0$ (we refer to the appendix for more details in the computations):
\begin{equation}
\begin{cases}
\begin{aligned}
&\p_{t}\rho-2\D\mu(\rho)+{\rm div}(\rho v)=0,\\
&\p_t {\rm div}v+ u\cdot\n {\rm div}v+\n v: ^{t}\n u-\frac{1}{2}\frac{\lambda(\rho)}{\rho^2}\n\rho\cdot{\rm div}{\rm curl}v\\
&\hspace{4cm}-R(\rho,v)-\frac{1}{2}\n\n\va(\rho):{\rm curl}v +\D F(\rho)=0,\\[2mm]
&\p_t{\rm curl}v+(u-\frac{1}{2}\n\va(\rho))\cdot\n{\rm curl}v-\frac{\mu(\rho)}{\rho}\D{\rm curl}v+R_1=0,\\[2mm]
&(\rho,{\rm div}v,{\rm curl}v)_{/t=0}=(\rho_{0},{\rm div}v_{0},{\rm curl}v_0).
\end{aligned}
\end{cases}
\label{1.b}
\end{equation}
with:
$$
\begin{aligned}
&R(\rho,v)=-\frac{1}{2}\p_i\va(\rho)\p_j({\rm curl}v)_{ij},\\
&(R^i_1)_{ij}=\sum_k (\p_{i} u_k\p_k v_j-\p_j u_k\p_k v_{i})-\frac{1}{2}\frac{\lambda(\rho)}{\rho^2}\big(\sum_k (\p_{i}\rho\p_k({\rm curl}v_{kj})-\p_j\rho\p_k({\rm curl }v_{ki})\big)\\
&\hspace{7cm}-\frac{1}{2}\sum_k \big(\p_{ik}\va(\rho)({\rm curl}v)_{kj}-\p_{kj}\va(\rho)({\rm curl}v)_{ik}\big).
\end{aligned}
$$
We can observe that $\rho$ and ${\rm curl}v$ have a parabolic behavior whereas 
${\rm div}v$ verifies a transport equation. In the sequel we are interested in proving existence of strong solution for the system (\ref{1.b}), let us start by recalling some results in the case of constant viscosity coefficients. Existence and uniqueness for (\ref{0.1}) in the case of smooth data with no vacuum has been stated in the pioneering works by Nash (see \cite{Nash}) and Matsumura and Nishida ( see \cite{MN1,MN}). These results have been extend in the framework of critical strong solution for the scaling of the equations  (see \cite{CD,DL,DW, arma,H1, H2} ). In the previous paper the authors use in a crucial way the parabolicity of the momentum equation on the velocity $u$ in order to get a gain of two derivatives in suitable Besov space; it allows to control the Lipschitz norm of $u$ (more precisely to show that $\n u$ belongs to $L^1(L^\infty)$) in order to estimate the density $\rho$ via the mass equation which verifies a transport equation. In particular it implies that in these works the density has no regularizing effects. In the case of the existence of global strong solution with small initial , it is proven in \cite{CD, arma} that the density has only a damping effect. More precisely if $\rho_0-1$ belongs to the Besov space $B^{\N}_{2,1}$ then the solution $\rho-1$ belongs in $L^1(\R^+,B^{\N}_{2,1})$.\\
Very few articles really take into account the structure of the viscosity coefficients in the framework of strong solution (we refer to \cite{BD,MV,CPAM,CPAM1} for the existence of global weak solution with viscosity coefficients verifying the relation (\ref{1.2}) or more precisely to the stability of the global weak solutions, indeed the problem of the existence of global weak solutions remains open). Let us mention however that in \cite{CH}  it is proven the existence of strong solution in critical Besov space, with assumption on the initial data which weaken the results of \cite{CD,DL,DW, arma,H1, H2}. Indeed $u_0$ is only in $B^{\N-1}_{2,2}\cap B^{-1}_{\infty,1}$. We would like also to mention the paper \cite{arxiv} where we obtain the existence of global strong solution for the shallow water system with large initial data for the initial data in terms of the scaling of the equation. To do this we consider the notion of quasi-solution (developed in the framework of global weak solutions in \cite{cras,CPAM,CPAM1}), which allows us in particular to deal with initial data $(\rho_0-1,u_0)$ large in $(B^{\N}_{2,\infty}\cap L^\infty)\times B^{\N-1}_{2,\infty}$ improving the smallness assumption of the previous paper. To do this we need to chose the initial velocity irrotational, it enables us to obtain global strong solution with large energy initial data in dimension $N=2$.\\
Now we are going to study the system (\ref{1.b}) and to show how obtain strong solution with large initial data in finite time and global strong solution with small initial data. The interest of this paper is not to weaken the regularity hypothesis on the initial data but  to prove that solution exists for system (\ref{0.1}) such that the density $\rho$ has regularizing effects. It is a priori very surprising since the system (\ref{0.1}) seems purely hyperbolic for the density, our result is really due to the specific structure of the viscosity coefficients which makes parabolic the system on the density. In other word there is a purely non linear effect which involves a parabolic behavior on the density. Let us mention that this choice of viscosity coefficients is crucial in \cite{BD,MV} in order to have additional regularity on the density (roughly speaking $\n\rho$ is in $L^\infty(L^2)$) which allows to have enough compactness to deal with the pressure and the viscosity coefficients.Let us observe that in counterpart we have a loss of regularity on ${\rm div}v$ compared with ${\rm div}u$since it behaves only as a transport equation (indeed if we consider $u$, then ${\rm div}u$ has a parabolic effects).\\
Let us briefly describe the ideas of the proof,  since ${\rm div}v$ behaves as a transport equation it is natural to choose an initial data set on ${\rm div}v_0$ which is embedded in $L^\infty$. Indeed in order to estimate ${\rm div}v$ via the second equation of (\ref{1.b}) which is a transport equation, we need to control $\n u$ in $L^1(L^\infty)$ in order to propagate the regularity of ${\rm div}v_0$ along the time.
To do this we take ${\rm div}v_0$ in $B^{\NN}_{p,1}$, now since we have $\D u={\rm div}{\rm curl}v+\n{\rm div}v-\n\D\va(\rho)$ it is important to control $\D\va(\rho)$ in $L^1(L^\infty)$. Using the regularizing effects of the heat equation it is sufficient to take $q_0$ in $B^{\NN}_{p,1}$. Considering the equation on ${\rm div}v$ and the coupling between $q$ and ${\rm curl}v$, it is natural to choose 
${\rm curl}v_0$ in $B^{\NN}_{p,1}$.
\\
To simplify the notation, we assume that $q=\rho-1$, hence as long as $\rho$ does not vanish, the equations for ($q$,$v$) read:
\begin{equation}
\begin{cases}
\begin{aligned}
&\p_t q-2\mu'(1+q)\D q-\mu''(1+q)|\n q|^2+{\rm div}(qv)+{\rm div}v=0,\\
&\p_t {\rm div}v+ u\cdot\n {\rm div}v+\n v: ^{t}\n u-\frac{1}{2}\frac{\lambda(1+q)}{(1+q)^2}\n\rho\cdot{\rm div}{\rm curl}v\\
&\hspace{4cm}-R(\va,v)-\frac{1}{2}\n\n\va(1+q):{\rm curl}v +\D F(1+q)=0,\\[2mm]
&\p_t{\rm curl}v+(u-\frac{1}{2}\n\va(1+q))\cdot\n{\rm curl}v-\frac{\mu(1+q)}{1+q}\D{\rm curl}v+R_1=0
\end{aligned}
\end{cases}
\label{0.6}
\end{equation}
In the sequel we will deal with the case $\mu(\rho)=\mu\rho^\alpha$ with $\alpha>1-\frac{1}{N}$ in order to verify the Lam\'e condition.
We can now state our main result:
\begin{theorem}
Let $N\geq 2$ and $1\leq p<+\infty$. Let $v_0\in L^\infty\cap B^{\NN+1}_{p,1}$, $q_0\in B^{\NN}_{p,1}$ and $\rho_0\geq c>0$ then it exists a time $T>0$ such that the system (\ref{0.6}) has a  solution $(q,v)$ with:
$$
\begin{cases}
\begin{aligned}
&v\in L^{\infty}_{T}(L^\infty),\;{\rm div} v\in \widetilde{C}_T(B^{\NN}_{p,1}),\\
&{\rm curl}v\in \widetilde{C}_{T}(B^{\NN}_{p,1})\cap \widetilde{L}_{T}^{1}(B^{\NN+2}_{p,1}) ,\\
&q\in \widetilde{C}_{T}(B^{\NN}_{p,1})\cap \widetilde{L}_{T}^{1}(B^{\NN+2}_{p,1}),\\
&(\rho,\frac{1}{\rho})\in L^{\infty}_{T}(L^{\infty}) .
\end{aligned}
\end{cases}
$$
In addition if $p<2N$ then the solution is unique.
\label{theo1}
\end{theorem}
\begin{remarka}
By definition of $v$ it implies that $u_0$ is in $B^{\NN-1}_{p,1}+B^{\NN+1}_{p,1}$ which is of course not optimal in terms of the scaling of the equations. In particular in terms of regularity on the initial data, this result is not optimal as in \cite{H2}. However we exhibit in a very surprising way some regularizing effects on the density which has a gain of two derivatives. This phenomena is really non linear and due to the specific choice on the viscosity coefficients.
\end{remarka}
\begin{remarka}
Let us however explain in which case the previous theorem is optimal in terms of regularity on the initial data. When we deal with the particular case $v_0=0$ which implies that $u_0=-\n\va(\rho_0)$  then in this situation $u_0$ is only in $B^{\NN-1}_{p,1}$ for any $1\leq p<+\infty$. It improves of two manners the existing results in \cite{H2,M3AS,DW}, the first thing is that we deal with larger space on the initial data since $u_0$ and $q_0$ are in $B^{\NN}_{p,1}$ and $B^{\NN-1}_{p,1}$ without any restriction on $p$. The second point is that we have regularizing effect on the density what is not the case in \cite{H2,M3AS,DW}.
\end{remarka}
Let us extend now the previous theorem to the case of global strong solution with small initial data.
\begin{theorem}
Let $N\geq 2$ and $1\leq p<\max(N,4)$. Let $v_0\in \widetilde{B}^{\N-1,\NN+1}_{2,p,1}$, $q_0\in  \widetilde{B}^{\N-1,\NN}_{2,p,1} $  then it exists $\e_{0}>0$ such that if:
\begin{equation}
\|q_0\|_{ \widetilde{B}^{\N-1,\NN}_{2,p,1} }+\|v_0\|_ { \widetilde{B}^{\N-1,\NN-1}_{2,p,1}}\leq \e_0,
\label{petitesse}
\end{equation}
then it exists a global strong solution $(q,v)$ of the system (\ref{0.6}) such that:
$$
\begin{cases}
\begin{aligned}
&{\rm div} v\in \widetilde{C}(\R^+,\widetilde{B}^{\N-2,\NN}_{2,p,1})\cap \widetilde{L}^1(\R^+, \widetilde{B}^{\NN,\NN}_{2,p,1})  ,\\
&{\rm curl}v\in \widetilde{C}(\R^+,\widetilde{B}^{\NN-2,\NN}_{2,p,1})\cap \widetilde{L}^{1}(\R^+,\widetilde{B}^{\N,\NN+2}_{2,p,1}) ,\\
&q\in \widetilde{C}(\R^+,\widetilde{B}^{\N-1,\NN}_{2,p,1})\cap \widetilde{L}^{1}(\R^+,\widetilde{B}^{\NN+1,\NN+2}_{2,p,1}),\\
&(\rho,\frac{1}{\rho})\in L^{\infty}_{T}(L^{\infty}) .
\end{aligned}
\end{cases}
$$
\label{theo2}
\end{theorem}


\begin{remarka}
We refer to \cite{arma} for the definition of the hybrid Besov spaces.
\end{remarka}
\begin{remarka}
Let us mention that the condition (\ref{petitesse}) is the same than in \cite{arma} and has a sense, indeed if $v_0\in\widetilde{B}^{\N-1,\NN+1}_{2,p,1}$ and $q_0\in \widetilde{B}^{\N-1,\NN}_{2,p,1}$ then $u_0$ is in $\widetilde{B}^{\N-1,\NN+1}_{2,p,1}$ which is embedded in $\widetilde{B}^{\N-1,\NN-1}_{2,p,1}$. The condition $p<\max(4,N)$ is technical and is due to the paraproduct law when we have Besov space with different behavior in low and high frequencies.
\end{remarka}
\begin{remarka}
Let us mention that this theorem is an extension of the result in \cite{arma,arxiv}. Indeed we have the same smallness assumption than in \cite{arma,arxiv}, and with additional regularity on the initial velocity since $u_0$ is in $B^{\NN+1}_{p,1}$ in high frequencies. However we are able to prove in this case a regularizing effects on the density which behaves as the solution of parabolic equation.
\end{remarka}
\begin{remarka}
In the particular case where $v_0=0$, it means $u_0=-\n\va(\rho_0)$ then our theorem improves the results of \cite{arma,CD} since we have the same initial data assumption but we prove in addition a regularizing effects on the density.
\end{remarka}

Our paper is structured as follows. In section \ref{section2}, we give a few notations and briefly introduce the basic Fourier analysis
techniques needed to prove our result. In section \ref{section3} and section \ref{section4}, we prove theorem \ref{theo1}  and more particular the existence of such solution in section \ref{section3} and the uniqueness in section \ref{section4}. In section \ref{section5} we are proving the global well-posedness of theorem \ref{theo1}.  We postpone in the appendix the proof of the equivalence between the system (\ref{0.1}) and the system (\ref{0.1a}).

\section{Littlewood-Paley theory and Besov spaces}
\label{section2}
As usual, the Fourier transform of $u$ with respect to the space variable will be denoted by $\mathcal{F}(u)$ or $\hat{u}$. 
In this section we will state classical definitions and properties concerning the homogeneous dyadic decomposition with respect to the Fourier variable. We will recall some classical results and we refer to \cite{BCD} (Chapter 2) for proofs (and more general properties).

To build the Littlewood-Paley decomposition, we need to fix a smooth radial function $\chi$ supported in (for example) the ball $B(0,\frac{4}{3})$, equal to 1 in a neighborhood of $B(0,\frac{3}{4})$ and such that $r\mapsto \chi(r.e_r)$ is nonincreasing over $\R_+$. So that if we define $\varphi(\xi)=\chi(\xi/2)-\chi(\xi)$, then $\varphi$ is compactly supported in the annulus $\{\xi\in \R^d, \frac{3}{4}\leq |\xi|\leq \frac{8}{3}\}$ and we have that,
\begin{equation}
 \forall \xi\in \R^d\setminus\{0\}, \quad \sum_{l\in\Z} \varphi(2^{-l}\xi)=1.
\label{LPxi}
\end{equation}
Then we can define the \textit{dyadic blocks} $(\ddl)_{l\in \Z}$ by $\ddl:= \varphi(2^{-l}D)$ (that is $\hat{\ddl u}=\varphi(2^{-l}\xi)\hat{u}(\xi)$) so that, formally, we have
\begin{equation}
u=\Sum_l \ddl u
\label{LPsomme} 
\end{equation}
As (\ref{LPxi}) is satisfied for $\xi\neq 0$, the previous formal equality holds true for tempered distributions \textit{modulo polynomials}. A way to avoid working modulo polynomials is to consider the set $\cS_h'$ of tempered distributions $u$ such that
$$
\lim_{l\rightarrow -\infty} \|\dot{S}_l u\|_{L^\infty}=0,
$$
where $\dot{S}_l$ stands for the low frequency cut-off defined by $\dot{S}_l:= \chi(2^{-l}D)$. If $u\in \cS_h'$, (\ref{LPsomme}) is true and we can write that $\dot{S}_l u=\Sum_{q\leq l-1} \ddq u$. We can now define the homogeneous Besov spaces used in this article:
\begin{defi}
\label{LPbesov}
 For $s\in\R$ and  
$1\leq p,r\leq\infty,$ we set
$$
\|u\|_{B^s_{p,r}}:=\bigg(\sum_{l} 2^{rls}
\|\Delta_l  u\|^r_{L^p}\bigg)^{\frac{1}{r}}\ \text{ if }\ r<\infty
\quad\text{and}\quad
\|u\|_{ B^s_{p,\infty}}:=\sup_{l} 2^{ls}
\|\Delta_l  u\|_{L^p}.
$$
We then define the space $ B^s_{p,r}$ as the subset of  distributions $u\in {\cS}'_h$ such that $\|u\|_{ B^s_{p,r}}$ is finite.
\end{defi}
Once more, we refer to \cite{BCD} (chapter $2$) for properties of the inhomogeneous and homogeneous Besov spaces. Among these properties, let us mention:
\begin{itemize}
\item for any $p\in[1,\infty]$ we have the following chain of continuous embeddings:
$$
 B^0_{p,1}\hookrightarrow L^p\hookrightarrow  B^0_{p,\infty};
$$
\item if $p<\infty$ then 
  $B^{\frac dp}_{p,1}$ is an algebra continuously embedded in the set of continuous 
  functions decaying to $0$ at infinity;
    \item for any  smooth homogeneous  of degree $m$ function $F$ on $\R^d\setminus\{0\}$
the operator $F(D)$ defined by $F(D)u=\mathcal{F}^{-1}\Big(F(\cdot)\mathcal{F}(u)(\cdot)\Big)$ maps  $ B^s_{p,r}$ in $ B^{s-m}_{p,r}.$ This implies that the gradient operator maps $ B^s_{p,r}$ in $B^{s-1}_{p,r}.$  
  \end{itemize}
We refer to \cite{BCD} (lemma 2.1) for the Bernstein lemma (describing how derivatives act on spectrally localized functions), that entails the following embedding result:
\begin{prop}\label{LP:embed}
\sl{For all $s\in\R,$ $1\leq p_1\leq p_2\leq\infty$ and $1\leq r_1\leq r_2\leq\infty,$
  the space $ B^{s}_{p_1,r_1}$ is continuously embedded in 
  the space $ B^{s-d(\frac1{p_1}-\frac1{p_2})}_{p_2,r_2}.$}
\end{prop}
Then we have:
$$
 B^\fdp_{p,1}\hookrightarrow  B^0_{\infty,1}\hookrightarrow L^\infty.
$$
In this paper, we shall mainly work with functions or distributions depending on both the time variable $t$ and the space variable $x.$ We shall denote by $\cC(I;X)$ the set of continuous functions on $I$ with values in $X.$ For $p\in[1,\infty]$, the notation $L^p(I;X)$ stands for the set of measurable functions on  $I$ with values in $X$ such that $t\mapsto \|f(t)\|_X$ belongs to $L^p(I)$.

In the case where $I=[0,T],$  the space $L^p([0,T];X)$ (resp. $\cC([0,T];X)$) will also be denoted by $L_T^p X$ (resp. $\cC_T X$). Finally, if $I=\R^+$ we shall alternately use the notation $L^p X.$

The Littlewood-Paley decomposition enables us to work with spectrally localized (hence smooth) functions rather than with rough objects. We naturally obtain bounds for each dyadic block in spaces of type $L^\rho_T L^p.$  Going from those type of bounds to estimates in  $L^\rho_T \dot B^s_{p,r}$ requires to perform a summation in $\ell^r(\Z).$ When doing so however, we \emph{do not} bound the $L^\rho_T \dot B^s_{p,r}$ norm for the time integration has been performed \emph{before} the $\ell^r$ summation.
This leads to the following notation:

\begin{defi}\label{d:espacestilde}
For $T>0,$ $s\in\R$ and  $1\leq r,\sigma\leq\infty,$
 we set
$$
\|u\|_{\tilde L_T^\sigma  B^s_{p,r}}:=
\bigl\Vert2^{js}\|\ddq u\|_{L_T^\sigma L^p}\bigr\Vert_{\ell^r(\Z)}.
$$
\end{defi}
One can then define the space $\tilde L^\sigma_T \dot B^s_{p,r}$ as the set of  tempered distributions $u$ over $(0,T)\times \R^d$ such that $\lim_{q\rightarrow-\infty}\dot S_q u=0$ in $L^\sigma([0,T];L^\infty(\R^d))$ and $\|u\|_{\tilde L_T^\sigma  B^s_{p,r}}<\infty.$ The letter $T$ is omitted for functions defined over $\R^+.$ 
The spaces $\tilde L^\sigma_T  B^s_{p,r}$ may be compared with the spaces  $L_T^\sigma \dot B^s_{p,r}$ through the Minkowski inequality: we have
$$
\|u\|_{\tilde L_T^\sigma  B^s_{p,r}}
\leq\|u\|_{L_T^\sigma  B^s_{p,r}}\ \text{ if }\ r\geq\sigma\quad\hbox{and}\quad
\|u\|_{\tilde L_T^\sigma  B^s_{p,r}}\geq
\|u\|_{L_T^\sigma  B^s_{p,r}}\ \text{ if }\ r\leq\sigma.
$$
All the properties of continuity for the product and composition which are true in Besov spaces remain true in the above  spaces. The time exponent just behaves according to H\"older's inequality. 
\medbreak
Let us now recall a few nonlinear estimates in Besov spaces. Formally, any product of two distributions $u$ and $v$ may be decomposed into 
\begin{equation}\label{eq:bony}
uv=T_uv+T_vu+R(u,v), \mbox{ where}
\end{equation}
$$
T_uv:=\sum_l\dot S_{l-1}u\ddl v,\quad
T_vu:=\sum_l \dot S_{l-1}v\ddl u\ \hbox{ and }\ 
R(u,v):=\sum_l\sum_{|l'-l|\leq1}\ddl u\,\dot\Delta_{l'}v.
$$
The above operator $T$ is called ``paraproduct'' whereas $R$ is called ``remainder''. The decomposition \eqref{eq:bony} has been introduced by Bony in \cite{BJM}.

In this article we will frequently use the following estimates (we refer to \cite{BCD} section 2.6):.
\begin{proposition}
Under the same assumptions there exists a constant $C>0$ such that if $1/p_1+1/p_2=1/p$, and $1/r_1+1/r_2=1/r$:
$$
\|\dot{T}_u v\|_{B_{2,1}^s}\leq C \|u\|_{L^\infty} \|v\|_{B_{2,1}^s},
$$
$$\|\dot{T}_u v\|_{B_{p,r}^{s+t}}\leq C\|u\|_{B_{p_1,r_1}^t} \|v\|_{\dot{B}_{p_2,r_2}^s} \quad (t<0),
$$
\begin{equation}
 \|\dot{R}(u,v)\|_{B_{p,r}^{s_1+s_2-\fd}} \leq C\|u\|_{B_{p_1,r_1}^{s_1}} \|v\|_{B_{p_2,r_2}^{s_2}} \quad (s_1+s_2>0).
\end{equation}
\label{produit}
\end{proposition}
Let us now turn to the composition estimates. We refer for example to \cite{BCD} (Theorem $2.59$, corollary $2.63$)):
\begin{proposition}
\sl{\begin{enumerate}
 \item Let $s>0$, $u\in B_{p,1}^s\cap L^{\infty}$ and $F\in W_{loc}^{[s]+2, \infty}(\R^d)$ such that $F(0)=0$. Then $F(u)\in B_{p,1}^s$ and there exists a function of one variable $C_0$ only depending on $s$, $p$, $d$ and $F$ such that
$$
\|F(u)\|_{B_{p,1}^s}\leq C_0(\|u\|_{L^\infty})\|u\|_{B_{p,1}^s}.
$$
\item If $u$ and $v\in B_{p,1}^\fd$ and if $v-u\in B_{p,1}^s$ for $s\in]-\fd, \fd]$ and $G\in W_{loc}^{[s]+3, \infty}(\R^d)$, then $G(v)-G(u)$ belongs to $B_{p,1}^s$ and there exists a function of two variables $C$ only depending on $s$, $d$ and $G$ such that
$$
\|G(v)-G(u)\|_{B_{p,1}^s}\leq C(\|u\|_{L^\infty}, \|v\|_{L^\infty})\left(|G'(0)| +\|u\|_{B_{p,1}^\fd} +\|v\|_{B_{p,1}^\fd}\right) \|v-u\|_{B_{p,1}^s}.
$$
\end{enumerate}}
\label{estimcompo}
\end{proposition}
Let us now recall a result of interpolation which explains the link between the space $B^{s}_{p,1}$ and the space $B^{s}_{p,\infty}$ (see
\cite{BCD} sections $2.11$ and $10.2.4$):
\begin{proposition}\sl{
\label{interpolationlog}
There exists a constant $C$ such that for all $s\in\R$, $\e>0$ and
$1\leq p<+\infty$,
$$\|u\|_{\widetilde{L}_{T}^{\sigma}(B^{s}_{p,1})}\leq C\frac{1+\e}{\e}\|u\|_{\widetilde{L}_{T}^{\sigma}(B^{s}_{p,\infty})}
\log\biggl(e+\frac{\|u\|_{\widetilde{L}_{T}^{\sigma}(B^{s-\e}_{p,\infty})}+ \|u\|_{\widetilde{L}_{T}^{\sigma}(B^{s+\e}_{p,\infty})}}
{\|u\|_{\widetilde{L}_{T}^{\sigma}(B^{s}_{p,\infty})}}\biggl).$$ \label{5Yudov}}
\end{proposition}
\subsection{Parabolic equations}
Let us end this section by recalling the following estimates for the heat equation:
\begin{proposition}\sl{
\label{chaleur} Let $s\in\R$, $(p,r)\in[1,+\infty]^{2}$ and
$1\leq\rho_{2}\leq\rho_{1}\leq+\infty$. Assume that $u_{0}\in B^{s}_{p,r}$ and $f\in\widetilde{L}^{\rho_{2}}_{T}
(B^{s-2+2/\rho_{2}}_{p,r})$.
Let u be a solution of:
$$
\begin{cases}
\begin{aligned}
&\p_{t}u-\mathcal{A} u=f\\
&u_{/t=0}=u_{0},\\
\end{aligned}
\end{cases}
$$  
with $\mathcal{A} u=-\mu\D u-(\lambda+\mu)\n{\rm div}u$ (where $\mu>0$ and $2\mu+\lambda>0$). Then there exists $C>0$ depending only on $N,\mu,\rho_{1}$ and
$\rho_{2}$ such that:
$$\|u\|_{\widetilde{L}^{\rho_{1}}_{T}(B^{s+2/\rho_{1}}_{p,r})}\leq C\big(
 \|u_{0}\|_{B^{s}_{p,r}}+\|f\|_{\widetilde{L}^{\rho_{2}}_{T}
 (B^{s-2+2/\rho_{2}}_{p})}\big)\,.$$
 If in addition $r$ is finite then $u$ belongs to $C([0,T],B^{s}_{p,r})$.}
\end{proposition}

Let us now state the following transport-diffusion estimates which are adaptations of the ones given in \cite{BCD} section 3, see also \cite{M3AS} for a proof. We consider here the following system:
\begin{equation}
\begin{cases}
\begin{aligned}
&\p_t v+u\cdot\n v-\mu b\D v=f,\\
&v(0,\cdot)=v_0,
\end{aligned}
\end{cases}
\label{E}
\end{equation}
with $b=1+a$ such that $a\in\widetilde{L}_T^\infty(B^{\NN}_{p,1})\cap \widetilde{L}_T^1(B^{\NN+2}_{p,1})$.
\begin{proposition}\sl{
 Let $1\leq p<+\infty$, $T>0$, $-\NN<s\leq \NN$, $v_0\in B_{p,1}^s$, $f\in L_T^1 B_{p,1}^s$ and $u\in \widetilde{L}_T^1(B^{\NN+1}_{p,1})$. If $v$ is a solution of (\ref{E}):
then setting $V(t)=\int_0^t( \|u(\tau)\|_{B^{\NN+1}_{p,1}}+\|a(\tau)\|^2_{B^{\NN+1}_{p,1}} )  d\tau$, there exists a constant $C>0$ such that for all $t\in[0,T]$:
$$
\|v\|_{\widetilde{L}_t^1(B^{s+2}_{p,1})} +\|v\|_{\widetilde{L}_t^\infty(B^{s}_{p,1})} \leq C e^{CV(t)} \Big(\|u_0\|_{\dot{B}_{p,1}^s}+ \int_0^t\|f(\tau)\|_{\dot{B}_{p,1}^s} e^{-CV(\tau)} d\tau\Big).
$$}
\label{lemmeuti}
\end{proposition}
Let us consider now the following transport diffusion equation:
\begin{equation}
\begin{cases}
\begin{aligned}
&\p_t q+v\cdot\n q-\mu \D q=-(1+q){\rm div}v,\\
&q(0,\cdot)=q_0,
\end{aligned}
\end{cases}
\label{E'}
\end{equation}

\begin{proposition}
Assume that $q_0\in B^{\NN}_{p,1}$, $v\in\widetilde{L}_T^1(B^{\NN+1}_{p,1})$ for $1\leq p<+\infty$ and $q\in \widetilde{L}_T^\infty(B^{\NN}_{p,1})$ satisfies (\ref{E'}). Let $V(t)=\int_0^t \|v(\tau)\|_{B^{\NN+1}_{p,1}}d\tau$, there exists a constant $C$ depending only on $N$ such that for all $t\in[0,T]$:
\begin{equation}
\|(Id-\dot{S}_m)q\|_{\widetilde{L}_T^\infty(B^{\NN}_{p,1})}\leq \|(Id-\dot{S}_m)q_0\|_{B^{\NN}_{p,1}}+
(e^{C V(t)}-1)(1+\| q_0\|_{B^{\NN}_{p,1}})
\label{utiltran}
\end{equation}
\label{chalHF}
\end{proposition}
\textbf{Proof: }  Localizing in frequency, if for $j\in \Z$, we have:
$$
\p_{t}\D_j u+v\cdot\n\D_j u-\mu\D\D_j u=R_j-\D_j((1+q){\rm div}v),
$$
where $R_j=[v.\n, \D_j]u$. Multiplying the equation by $|\D_j q|^{p-2}\D_j q$, integrating by parts, using H\"older's inequality and lemma A5 in \cite{DL}, we obtain with $\alpha>0$ depending on $\mu$ and $p$:
$$
\frac{1}{p}\p_{t}\|\D_j u|^p_{L^p}+\alpha 2^{2j} \|u_j\|_{L^p}^p \leq \|\D_j u\|_{L^p}^{p-1}(\|R_j\|_{L^p}+\|\D_j((1+q){\rm div}v)\|_{L^p}+\frac{1}{p}\|\D_j u\|_{L^p}\|{\rm div}v\|_{L^\infty})).
$$
Times integration leads to:
$$
\begin{aligned}
&\|\D_j u(t)\|_{L^p}+\alpha 2^{2j} \int^t_0 \|u_j(\tau)\|_{L^p}d\tau \\
&\leq \|\D_j q_0\|_{L^p}+\int^t_0\big(\|R_j\|_{L^p}+\|\D_j((1+q){\rm div}v)\|_{L^p}+\frac{1}{p}\|\D_j u\|_{L^p}\|{\rm div}v\|_{L^\infty}\big)d\tau.
\end{aligned}
$$
Classical commutator estimates (we refer to \cite{BCD} section $2.10$) then imply that there exists a summable positive sequence $c_j=c_j(t)$ whose sum is $1$ such that:
$$
\|R_j\|_{L^p}\leq c_j 2^{-js} \|\n v\|_{B^{\NN}_{p,\infty}\cap L^\infty} \|u\|_{B^{\NN}_{p,1}}.
$$
Thanks to the paraproduct and remainder laws, we have:
$$
\|\D_j((1+q){\rm div}v\|_{L^p}\lesssim C c_{j}(t)2^{-j\NN}(1+\|q\|_{B^\NN_{p,1}})\|{\rm div}v\|_{B^{\NN}_{p,1}}.
$$
We get finally for all $t\in[0,T]$:
$$
\begin{aligned}
&2^{j\NN}\|\D_j q(t)\|_{L^p}+\alpha 2^{j(2+\NN)} \int^t_0 \|u_j(\tau)\|_{L^p}d\tau\leq2^{j\NN}\|\D_j q_0\|_{L^p}+C\int^t_0 c_j(1+\|q\|_{B^{\NN}_{p,1}}) V' d\tau.
\end{aligned}
$$
Summing up on $j\geq m$ and using the fact that (to see this it suffices to apply Gronwall inequality to the previous inequality):
$$\|q\|_{\widetilde{L}^\infty_t(B^{\NN}_{p,1})}\leq e^{C V(t)}\|q_0\|_{B^{\NN}_{p,1}}+e^{CV(t)}-1,$$
it gives:
$$
\begin{aligned}
\sum_{j\geq m}2^{j\NN}\|\D_j q(t)\|_{L^p}\leq \sum_{j\geq m}2^{j\NN}\|\D_j q_0\|_{L^p}
+\int^t_0 C V'(e^{C V}\|q_0\|_{B^{\NN}_{p,1}}+e^{CV}-1)d\tau +\int^t_0 CV' d\tau.
\end{aligned}
$$
Straightforward calculations lead to (\ref{utiltran}).
 $\blacksquare$

\subsection{Transport equations}
We begin this section by recalling some estimates in Besov spaces for transport and heat equations. For more details, the reader is referred to \cite{BCD}.
\begin{proposition}
Let $1\leq p_{1}\leq p\leq+\infty$, $r\in[1,+\infty]$ and $s\in\R$ be such that:
$$-N\min(\frac{1}{p_{1}},\frac{1}{p^{'}})<s<1+\frac{N}{p_{1}}.$$
Suppose that $q_{0}\in B^{s}_{p,r}$, $F\in L^{1}(0,T, B^{s}_{p,r})$ and that $q\in
L^{\infty}_{T}(B^{s}_{p,r})\cap
C([0,T];{\cal S}^{'})$ solves the following transport equation:
$$
\begin{cases}
\begin{aligned}
&\p_{t}q+u\cdot\n q=F,\\
&q_{\ t=0}=q_{0}.
\end{aligned}
\end{cases}
$$
There exists a constant $C$ depending only on $N$, $p$, $p_{1}$, $r$ and $s$ such that , we have for a.e $t\in[0,T]$:
\begin{equation}
\|q\|_{\widetilde{L}^{\infty}_{t}(B^{s}_{p,r})}\leq e^{CU(t)}\big(\|q_{0}\|_{B^{s}_{p,r}}+\int^{t}_{0}e^{-CU(\tau)}
\|F(\tau)\|_{B^{s}_{p,r}}d\tau\big),
\label{20}
\end{equation}
with:
$U(t)=\int^{t}_{0}\|\n u(\tau)\|_{B^{\frac{N}{p_{1}}}_{p_{1},\infty}\cap L^{\infty}}d\tau$.
\label{transport1}
\end{proposition}

\section{Existence of solution}
\label{section3}
\subsection{A priori estimates}
\label{aprioriqu}
Let us emphasize that we use the fact that $P'(1)>0$ only in the global existence result. Denoting by $G$ the unique primitive of $x\mapsto P'(x)/x$ such that $G(1)=0$, recall that system (\ref{0.6}) now reads:
\begin{equation}
\begin{cases}
\begin{aligned}
&\p_t q-2\mu'(1+q)\D q-2\mu''(1+q)|\n q|^2+{\rm div}(qv)+{\rm div}v=0,\\
&\p_t {\rm div}v+ u\cdot\n {\rm div}v+\n v: ^{t}\n u-\frac{1}{2}\frac{\lambda(1+q)}{(1+q)^2}\n q\cdot{\rm div}{\rm curl}v\\
&\hspace{4cm}-R(\rho,v)-\frac{1}{2}\n\n\va(1+q):{\rm curl}v +\D F(1+q)=0,\\[2mm]
&\p_t{\rm curl}v+(u-\frac{1}{2}\n\va(1+q))\cdot\n{\rm curl}v-\frac{\mu(1+q)}{1+q}\D{\rm curl}v+R_1=0
\end{aligned}
\end{cases}
\label{syst}
\end{equation}
with:
$$
\begin{aligned}
&R(\rho,v)=-\frac{1}{2}\p_i\va(\rho)\p_j({\rm curl}v)_{ij},\\
&(R_1)_{ij}=\sum_k (\p_i u_k\p_k v_j-\p_j u_k\p_k v_i)-\frac{1}{2}\frac{\lambda(\rho)}{\rho^2}\big(\sum_k (\p_i\rho\p_k({\rm curl}_{kj})-\p_j\rho\p_k({\rm curl}_{ki})\big)\\
&\hspace{7cm}-\frac{1}{2}\sum_k \big(\p_{ik}\va(\rho)({\rm curl}v)_{kj}-\p_{kj}\va(\rho)({\rm curl}v)_{ik}\big).
\end{aligned}
$$
In the sequel we are going to consider the case $\mu(\rho)=\mu\rho$ in order to simplify the calculus.
Let $(q_L,v_L)$ be the unique global solution which verifies ${\rm div}v_L=0$ and:
\begin{equation}
\begin{cases}
\begin{aligned}
&\p_{t}q_L-2\mu\D q_L=0,\\
&\p_{t}{\rm curl}v_L-\mu\D{\rm curl}v_L=0,\\
&(q_L,{\rm curl}v_L)_{/t=0}=(q_0,{\rm curl}v_0).
\end{aligned}
\end{cases}
\label{systlin}
\end{equation}
\begin{remarka}
In order to deal with the generic case $\mu(\rho)=\mu\rho$ it suffices to study the nonlinear system where $v_L$ verifies ${\rm div}v_L=0$ and:
\begin{equation}
\begin{cases}
\begin{aligned}
&\p_{t}q_L-2(\mu'(1)+\dot{S}_m(\mu'(1+q)-\mu'))\D q_L=0,\\
&\p_{t}{\rm curl}v_L-\mu\D{\rm curl}v_L=0,\\
&(q_L,{\rm curl}v_L)_{/t=0}=(q_0,{\rm curl}v_0).
\end{aligned}
\end{cases}
\label{systlin1}
\end{equation}
To do this we can use the proposition \ref{lemmeuti}.
\end{remarka}
Thanks to the classical heat estimates recalled in Proposition \ref{chaleur} (we refer for example to \cite{BCD}, lemma $2.4$ and chapter $3$), as $v_0\in B^{\NN}_{p,1}$ we have for all time $t$ and $C>0$:
\begin{equation}
\begin{aligned}
&\|{\rm curl}v_L\|_{\widetilde{L}_t^\infty(B^{\NN}_{p,1})}+\|q_L\|_{\widetilde{L}_t^\infty(B^{\NN}_{p,1})}\leq C (\|{\rm curl}v_0\|_{B^\NN_{p,1}}+\|q_0\|_{B^\NN_{p,1}})=CE_0.
\end{aligned}
\label{reguL1}
\end{equation}
Furthermore for a small $\eta>0$ (that we will precise in the sequel), it exists a time $T$ such that for $1\leq r<+\infty$:
\begin{equation}
\begin{aligned}
&\|{\rm curl}v_L\|_{\widetilde{L}^1_T (B^{\NN+2}_{p,1})}+\|q_L\|_{\widetilde{L}^1_T (B^{\NN+2}_{p,1})}  \leq \eta,\\
&\|{\rm curl}v_L\|_{\widetilde{L}_T^r (B^{\NN+\frac{2}{r}}_{p,1})}+  \|q_L\|_{\widetilde{L}_T^r (B^{\NN+\frac{2}{r}}_{p,1})}\leq C \eta^{\frac{1}{r}} E_0^{1-\frac{1}{r}}.
\end{aligned}
\label{reguL2}
\end{equation}
We are going now to set:
$${\rm curl}v={\rm curl}v_L+{\rm curl}\bar{v};\;{\rm div}\bar{v}={\rm div}v\;\;\mbox{and}\;\;q=q_L+\bar{q}.$$
We have to estimates now in Besov space the unknown $(\bar{q},{\rm div}\bar{v},{\rm curl}\bar{v})$, let us start with rewriting system (\ref{syst}) in terms of these unknown:
\begin{equation}
\begin{cases}
\begin{aligned}
&\p_t \bar{q}-2\mu \D \bar{q}
=F,\\
&\p_t {\rm div}\bar{v}+ u\cdot\n {\rm div}\bar{v}=G,\\
&\hspace{4cm}+R(\va,v)+\frac{1}{2}\n\n\va(1+q):{\rm curl}v -\D F(1+q),\\[2mm]
&\p_t {\rm curl}\bar{v}+u\cdot\n{\rm curl}v-\big(\mu+\dot{S}_m\big(\frac{\mu(1+q)}{1+q}-\mu(1)\big)\big)\D{\rm curl}\bar{v}=H.
\end{aligned}
\end{cases}
\label{syst1}
\end{equation}
with:
$$
\begin{aligned}
&F=-{\rm div}((q+1)v),\\
&G=-\n v: ^{t}\n u+\frac{1}{2}\frac{\lambda(1+q)}{(1+q)^2}\n q\cdot{\rm div}{\rm curl}v+R(\rho,v)+\frac{1}{2}\n\n\va(1+q):{\rm curl}v -\D F(1+q),\\
&H=\frac{1}{2}\n\va(1+q)\cdot\n{\rm curl}v+(Id-\dot{S}_m)\big(\frac{\mu(1+q)}{1+q}-\mu\big)\D{\rm curl}\bar{v}+(\frac{\mu(1+q)}{1+q}-\mu)\D{\rm curl}v_L-R_1.
\end{aligned}
$$
As we want to estimate on $(q,{\rm div}\bar{v},{\rm curl}\bar{v})$ in Besov spaces, we are going to use the propositions \ref{chaleur},  \ref{lemmeuti} applied to $q$ and ${\rm curl}\bar{v}$ and the proposition \ref{transport1} applied to ${\rm div}\bar{v}$. More precisely we have 
using proposition \ref{chaleur} and \ref{produit}:
\begin{equation}
\begin{aligned}
&\|\bar{q}\|_{\widetilde{L}_t^1(B^{\NN+2}_{p,1})}+\|\bar{q}\|_{\widetilde{L}_t^\infty (B^{\NN}_{p,1})}
\leq C\int^t_0 \| {\rm div}((q+1)v)(s,\cdot)\|_{B^\NN_{p,1}} ds.\\[2mm]
&C\leq \int^t_0\big(  \|q\|_{B^{\NN+1}_{p,1}}\|v\|_{L^\infty}+  \|v\|_{B^{\NN+1}_{p,1}}\|q\|_{L^\infty}+\|v\|_{B^{\NN+1}_{p,1}}  )  ds\\
&C\leq (\sqrt{t}\|q\|_{\widetilde{L}_t^2(B^{\NN+1}_{p,1})}\|v\|_{L_t^{\infty}(L^\infty)}+
\sqrt{t}\|v\|_{\widetilde{L}_t^2(B^{\NN+1}_{p,1})}\|q\|_{\widetilde{L}_t^\infty(B^{\NN}_{p,1})}+t\|v\|_{\widetilde{L}_t^\infty(B^{\NN+1}_{p,1})}).
\end{aligned}
\label{estimq2}
\end{equation}
It remains to estimate $\|v\|_{L_t^\infty(L^\infty)}$ by using the momentum equation, we observe that $v$ verifies a transport equation with a rest of the form $-\n F(\rho)+\frac{\mu}{\rho}{\rm div}(\rho{\rm curl}v)$. We get by using proposition \ref{estimcompo}:
\begin{equation}
\begin{aligned}
&\|v\|_{L_t^\infty(L^\infty)}\leq (\|v_0\|_{L^\infty}+\|\frac{1}{\rho}{\rm div}(\mu(\rho){\rm curl}v)\|_{\widetilde{L}_t^1(B^{\NN}_{p,1})}+\|\n F(\rho)\|_{\widetilde{L}_t^1(B^{\NN}_{p,1})}),\\
&\leq  \biggl(\|v_0\|_{L^\infty}+(1+C[\|q\|_{L_t^\infty}]\|q\|_{\widetilde{L}^\infty_t(B^{\NN}_{p,1})}) \big(\|{\rm curl} v\|_{\widetilde{L}_t^1(B^{\NN+1}_{p,1})}C[1+\|q\|_{L^\infty}]\\
&+\|{\rm curl}v\|_{\widetilde{L}_t^{\infty}(B^{\NN}_{p,1})}(1+ \|q\|_{\widetilde{L}_t^1(B^{\NN+1}_{p,1})}C[\|q\|_{L_t^\infty(L^\infty)}])\big)+C[\|q\|_{L_t^\infty(L^\infty)}]\|q\|_{\widetilde{L}_t^1(B^{\NN+1}_{p,1})}\biggl).
\end{aligned}
\label{controlev}
\end{equation}
Combining (\ref{estimq2}) and (\ref{controlev}) we obtain:
\begin{equation}
\begin{aligned}
&\|\bar{q}\|_{\widetilde{L}_t^1(B^{\NN+2}_{p,1})}+\|\bar{q}\|_{\widetilde{L}_t^\infty (B^{\NN}_{p,1})}\\
&\leq C \biggl(\sqrt{t}\|q\|_{\widetilde{L}^2(B^{\NN+1}_{p,1})}\biggl(\|v_0\|_{L^\infty}+(1+C[\|q\|_{L_t^\infty}]\|q\|_{\widetilde{L}^\infty_t(B^{\NN}_{p,1})}) \big(\|{\rm curl} v\|_{\widetilde{L}_t^1(B^{\NN+1}_{p,1})}C[1+\|q\|_{L^\infty}]\\
&+\|{\rm curl}v\|_{\widetilde{L}_t^{\infty}(B^{\NN}_{p,1})}(1+ \|q\|_{\widetilde{L}_t^1(B^{\NN+1}_{p,1})}C[\|q\|_{L_t^\infty(L^\infty)}])\big)+C[\|q\|_{L_t^\infty(L^\infty)}]\|q\|_{\widetilde{L}_t^1(B^{\NN+1}_{p,1})}\biggl) \\
&\hspace{6cm}+
\sqrt{t}\|v\|_{\widetilde{L}^2(B^{\NN+1}_{p,1})}\|q\|_{\widetilde{L}^\infty(B^{\NN}_{p,1})}+t\|v\|_{\widetilde{L}_t^\infty(B^{\NN+1}_{p,1})}\biggl).
\end{aligned}
\label{estimq3}
\end{equation}
Proceeding similarly and using proposition \ref{lemmeuti}, we have:
\begin{equation}
\begin{aligned}
\|{\rm curl}\bar{v}\|_{\widetilde{L}_t^\infty (B^{\NN}_{p,1})}+\|{\rm curl}\bar{v}\|_{\widetilde{L}_t^\infty (B^{\NN}_{p,1})}&\leq Ce^{C V_2(t)}\int^t_0 \|H(s)\|_{B^\NN_{p,1}}ds,
\end{aligned}
\label{estimv2}
\end{equation}
with $V_2(t)=\int^t_0(\|u(\tau)\|_{B^{\NN+1}_{p,1}}+\|q(\tau)\|^2_{B^{\NN+1}_{p,1}})d\tau$.
It suffices now to estimate $\|H(s)\|_{B^\NN_{p,1}}$ by using proposition \ref{produit} and \ref{estimcompo}. We have then:
\begin{equation}
\begin{aligned}
&\|\n\va(1+q)\cdot\n{\rm curl}v\|_{B^\NN_{p,1}}\lesssim C[\|q\|_{L^\infty}]\|q\|_{B^{\NN+1}_{p,1}}({\rm curl}\bar{v}\|_{B^{\NN+1}_{p,1}}+\|{\rm curl}v_L\|_{B^{\NN+1}_{p,1}}),\\
&\|(Id-\dot{S}_m)\big(\frac{\mu(1+q)}{1+q}-\mu\big)\D{\rm curl}\bar{v}\|_{B^\NN_{p,1}} \lesssim C[\|q\|_{L^\infty}]\|(Id-\dot{S}_m)q\|_{B^\NN_{p,1}}\|{\rm curl}\bar{v}\|_{B^{\NN+2}_{p,1}},\\
&\|(\frac{\mu(1+q)}{1+q}-\mu)\D{\rm curl}v_L\|_{B^\NN_{p,1}}\lesssim C[\|q\|_{L^\infty}]\|q\|_{B^\NN_{p,1}}\|{\rm curl}v_L\|_{B^{\NN+2}_{p,1}},\\
&\|R_1\|_{B^\NN_{p,1}}\lesssim  (\|{\rm div}\bar{v}\|_{B^\NN_{p,1}}+\|{\rm curl}\bar{v}\|_{B^\NN_{p,1}}+\|{\rm curl}v_L\|_{B^\NN_{p,1}})^2\\
&+ (1+C[\|q\|_{L^\infty}]\|q\|_{B^\NN_{p,1}})\|q\|_{B^{\NN+1}_{p,1}}(\|{\rm curl}\bar{v}\|_{B^{\NN+1}_{p,1}}
+ \|{\rm curl} v_L\|_{B^{\NN+1}_{p,1}})\\
&\hspace{4cm}+ C[\|q\|_{L^\infty}]\|q\|_{B^{\NN+2}_{p,1}}(\|{\rm curl}\bar{v}\|_{B^{\NN}_{p,1}}
+ \|{\rm curl} v_L\|_{B^{\NN}_{p,1}}).
\end{aligned}
\end{equation}
Finally we have obtained:
\begin{equation}
\begin{aligned}
&\|{\rm curl}\bar{v}\|_{\widetilde{L}_t^\infty (B^{\NN}_{p,1})}+\|{\rm curl}\bar{v}\|_{\widetilde{L}_t^1 (B^{\NN+2}_{p,1})}\leq Ce^{C V_2(t)}\big(\\
&+\int^t_0 C[\|q\|_{L^\infty}]\|q\|_{B^{\NN+1}_{p,1}}(\|{\rm curl}\bar{v}\|_{B^{\NN+1}_{p,1}}+\|{\rm curl}v_L\|_{B^{\NN+1}_{p,1}})  ds \\
&+ \int^t_0\big(C[\|q\|_{L^\infty}]\|(Id-\dot{S}_m)q\|_{B^\NN_{p,1}}\|{\rm curl}\bar{v}\|_{B^{\NN+2}_{p,1}}+  C[\|q\|_{L^\infty}]\|q\|_{B^\NN_{p,1}}\|{\rm curl}v_L\|_{B^{\NN+2}_{p,1}} \big)ds\\
&+\int^t_0  (\|{\rm div}\bar{v}\|_{B^\NN_{p,1}}+\|{\rm curl}\bar{v}\|_{B^\NN_{p,1}}+\|{\rm curl}v_L\|_{B^\NN_{p,1}})^2 ds\\
&+ \int^t_0C[\|q\|_{L^\infty}]\|q\|_{B^{\NN+2}_{p,1}}(\|{\rm curl}\bar{v}\|_{B^{\NN}_{p,1}}
+ \|{\rm curl} v_L\|_{B^{\NN}_{p,1}})ds
\\
&+\int^t_0 (1+C[\|q\|_{L^\infty}]\|q\|_{B^\NN_{p,1}})\|q\|_{B^{\NN+1}_{p,1}}(\|{\rm curl}\bar{v}\|_{B^{\NN+1}_{p,1}}
+ \|{\rm curl} v_L\|_{B^{\NN+1}_{p,1}})ds.
\end{aligned}
\label{estimcurl}
\end{equation}
We deduce that:
\begin{equation}
\begin{aligned}
&\|{\rm curl}\bar{v}\|_{\widetilde{L}_t^\infty (B^{\NN}_{p,1})}+\|{\rm curl}\bar{v}\|_{\widetilde{L}_t^1(B^{\NN+2}_{p,1})}\leq Ce^{C V_2(t)}\big(\\
& C[\|q\|_{L_t^\infty(L^\infty)}]\|q\|_{\widetilde{L}_t^2(B^{\NN+1}_{p,1})}(\|{\rm curl}\bar{v}\|_{\widetilde{L}_t^2(B^{\NN+1}_{p,1})}+\|{\rm curl}v_L\|_{\widetilde{L}_t^2(B^{\NN+1}_{p,1})})   \\
&+C[\|q\|_{L_t^\infty(L^\infty)}]\|(Id-\dot{S}_m)q\|_{\widetilde{L}_t^\infty(B^\NN_{p,1})}\|{\rm curl}\bar{v}\|_{\widetilde{L}_t^1(B^{\NN+2}_{p,1})}\\
&+  C[\|q\|_{L_t^\infty(L^\infty)}]\|q\|_{\widetilde{L}_t^{\infty}(B^\NN_{p,1})}\|{\rm curl}v_L\|_{\widetilde{L}_t^1(B^{\NN+2}_{p,1})}\\
&+(\sqrt{t}\|{\rm div}\bar{v}\|_{\widetilde{L}_t^{\infty}(B^\NN_{p,1})}+\sqrt{t}\|{\rm curl}\bar{v}\|_{\widetilde{L}_t^\infty (B^\NN_{p,1})}+\sqrt{t}\|{\rm curl}v_L\|_{\widetilde{L}_t^\infty(B^\NN_{p,1})})^2\\
&+C[\|q\|_{L_t^\infty(L^\infty)}]\|q\|_{\widetilde{L}_t^1(B^{\NN+2}_{p,1})}(\|{\rm curl}\bar{v}\|_{\widetilde{L}_t^\infty (B^{\NN}_{p,1})}
+ \|{\rm curl} v_L\|_{\widetilde{L}_t^\infty(B^{\NN}_{p,1})})
\\
&+ (1+C[\|q\|_{L_t^\infty(L^\infty)}]\|q\|_{\widetilde{L}_t^\infty (B^\NN_{p,1})})\|q\|_{\widetilde{L}_t^2(B^{\NN+1}_{p,1})}\\
&\hspace{4cm}\times(\|{\rm curl}\bar{v}\|_{\widetilde{L}_t^2(B^{\NN+1}_{p,1})}
+ \|{\rm curl} v_L\|_{\widetilde{L}_t^2(B^{\NN+1}_{p,1})}).
\end{aligned}
\label{estimcurla}
\end{equation}
We are going now to use the proposition \ref{utiltran} with $V(t)=\int^t_0\|v(\tau)\|_{B^{\NN+1}_{p,1}}d\tau$ which yields:
\begin{equation}
\begin{aligned}
&\|{\rm curl}\bar{v}\|_{\widetilde{L}_t^\infty (B^{\NN}_{p,1})}+\|{\rm curl}\bar{v}\|_{\widetilde{L}_t^1(B^{\NN+2}_{p,1})}\leq Ce^{C V_2(t)}\big(\\
& C[\|q\|_{L^\infty_t(L^\infty)})\|q\|_{\widetilde{L}_t^2(B^{\NN+1}_{p,1})}(\|{\rm curl}\bar{v}\|_{\widetilde{L}_t^2(B^{\NN+1}_{p,1})}+\|{\rm curl}v_L\|_{\widetilde{L}_t^2(B^{\NN+1}_{p,1})})   \\
&+C[\|q\|_{L_t^\infty(L^\infty)}]\big(\|(Id-S_m)q_0\|_{B^{\NN}_{p,1}}+(e^{C V(t)}-1)(1+\|q_0\|_{B^{\NN}_{p,1}})        \big)\|{\rm curl}\bar{v}\|_{\widetilde{L}_t^1(B^{\NN+2}_{p,1})}\\
&+  C[\|q\|_{L_t^\infty(L^\infty)}]\|q\|_{\widetilde{L}_t^{\infty}(B^\NN_{p,1})}\|{\rm curl}v_L\|_{\widetilde{L}_t^1(B^{\NN+2}_{p,1})}\\
&+(\sqrt{t}\|{\rm div}\bar{v}\|_{\widetilde{L}_t^{\infty}(B^\NN_{p,1})}+\sqrt{t}\|{\rm curl}\bar{v}\|_{\widetilde{L}_t^\infty (B^\NN_{p,1})}+\sqrt{t}\|{\rm curl}v_L\|_{\widetilde{L}_t^\infty(B^\NN_{p,1})})^2\\
&+C[\|q\|_{L_t^\infty(L^\infty)}]\|q\|_{\widetilde{L}_t^1(B^{\NN+2}_{p,1})}(\|{\rm curl}\bar{v}\|_{\widetilde{L}_t^\infty (B^{\NN}_{p,1})}
+ \|{\rm curl} v_L\|_{\widetilde{L}_t^\infty(B^{\NN}_{p,1})})
\\
&+ (1+C[\|q\|_{L_t^\infty(L^\infty)}]\|q\|_{\widetilde{L}_t^\infty (B^\NN_{p,1})})\|q\|_{\widetilde{L}_t^2(B^{\NN+1}_{p,1})}(\|{\rm curl}\bar{v}\|_{\widetilde{L}_t^2(B^{\NN+1}_{p,1})}
+ \|{\rm curl} v_L\|_{\widetilde{L}_t^2(B^{\NN+1}_{p,1})}).
\end{aligned}
\label{estimcurlab}
\end{equation}
In the same way, using proposition \ref{transport1} we have:
\begin{equation}
\begin{aligned}
\|{\rm div}\bar{v}\|_{\widetilde{L}_t^\infty (B^{\NN}_{p,1})}&\leq e^{C V_1(t)}(\|{\rm div}\bar{v}_0\|_{B^\NN_{p,1}}+\int^t_0 \|G(s)\|_{B^\NN_{p,1}}ds).
\end{aligned}
\label{estimv1}
\end{equation}
with $V_1(t)=\int^t_0\|u(s)\|_{B^{\NN+1}_{p,1}}ds$. It remains to estimate $ \|G(s)\|_{B^\NN_{p,1}}$ by using proposition \ref{produit} and \ref{estimcompo}, more precisely we have:
\begin{equation}
\begin{aligned}
&\|\n v: ^{t}\n u\|_{B^\NN_{p,1}}\lesssim (\|{\rm div}\bar{v}\|_{B^\NN_{p,1}}+\|{\rm curl}\bar{v}\|_{B^\NN_{p,1}}+\|{\rm curl}v_L\|_{B^\NN_{p,1}})(\|{\rm div}\bar{v}\|_{B^\NN_{p,1}}+\|{\rm curl}v\|_{B^\NN_{p,1}}+\|q\|_{B^{\NN+2}_{p,1}}),\\
&\|\frac{\lambda(1+q)}{(1+q)^2}\n q\cdot{\rm div}{\rm curl}v\|_{B^\NN_{p,1}}\lesssim 
(1+C[\|q\|_{L^\infty}]\|q\|_{B^\NN_{p,1}})\|q\|_{B^{\NN+1}_{p,1}}\\
&\hspace{5cm}\times(\|{\rm curl}\bar{v}\|_{B^{\NN+1}_{p,1}}
+ \|{\rm curl} v_L\|_{B^{\NN+1}_{p,1}}),\\
&\|R(\rho,v)\|_{B^\NN_{p,1}}\lesssim (1+C[\|q\|_{L^\infty}]\|q\|_{B^\NN_{p,1}})\|q\|_{B^{\NN+1}_{p,1}}(\|{\rm curl}\bar{v}\|_{B^{\NN+1}_{p,1}}
+ \|{\rm curl} v_L\|_{B^{\NN+1}_{p,1}})\\
 &\|\D F(1+q)\|_{B^\NN_{p,1}}\lesssim C[\|q\|_{L^\infty}]\|q\|_{B^{\NN+2}_{p,1}}.
\end{aligned}
\end{equation}
In particular it implies that:
\begin{equation}
\begin{aligned}
&\|{\rm div}\bar{v}\|_{\widetilde{L}_t^\infty (B^{\NN}_{p,1})}
\leq Ce^{C V_1(t)}\big(\|{\rm div}v_0\|_{B^\NN_{p,1}}+\int^t_0(1+C[\|q\|_{L^\infty}]\|q\|_{B^\NN_{p,1}})\|q\|_{B^{\NN+1}_{p,1}}\\
&\times(\|{\rm curl}\bar{v}\|_{B^{\NN+1}_{p,1}}
+ \|{\rm curl} v_L\|_{B^{\NN+1}_{p,1}})ds+\int^t_0 C[\|q\|_{L^\infty}]\|q\|_{B^{\NN+2}_{p,1}}ds\\
&+\int^t_0   (\|{\rm div}\bar{v}\|_{B^\NN_{p,1}}+\|{\rm curl}\bar{v}\|_{B^\NN_{p,1}}+\|{\rm curl}v_L\|_{B^\NN_{p,1}})(\|{\rm div}\bar{v}\|_{B^\NN_{p,1}}+\|{\rm curl}v\|_{B^\NN_{p,1}}+\|q\|_{B^{\NN+2}_{p,1}})ds\\
&\hspace{3cm}+\int^t_0  C[\|q\|_{L^\infty}]\|q\|_{B^{\NN+2}_{p,1}}(\|{\rm curl}\bar{v}\|_{B^{\NN}_{p,1}}+\|{\rm curl}v_L\|_{B^{\NN}_{p,1}})ds  \big),
\end{aligned}
\label{estimdiv}
\end{equation}
This gives:
\begin{equation}
\begin{aligned}
&\|{\rm div}\bar{v}\|_{\widetilde{L}_t^\infty (B^{\NN}_{p,1})}
\leq Ce^{C V_1(t)}\big(\|{\rm div}v_0\|_{B^\NN_{p,1}}+(1+C[\|q\|_{L_t^\infty(L^\infty)})\|q\|_{\widetilde{L}_t^\infty(B^\NN_{p,1})})\|q\|_{\widetilde{L}_t^2(B^{\NN+1}_{p,1})}\\
&\times(\|{\rm curl}\bar{v}\|_{\widetilde{L}_t^2(B^{\NN+1}_{p,1})}
+ \|{\rm curl} v_L\|_{\widetilde{L}^2(B^{\NN+1}_{p,1})})+C[\|q\|_{L_t^\infty(L^\infty)}]\|q\|_{\widetilde{L}_t^1(B^{\NN+2}_{p,1})}\\
&+ (\sqrt{t}\|{\rm div}\bar{v}\|_{\widetilde{L}_t^\infty(B^\NN_{p,1})}+\sqrt{t}\|{\rm curl}\bar{v}\|_{\widetilde{L}_t^\infty(B^\NN_{p,1})}+\sqrt{t}\|{\rm curl}v_L\|_{\widetilde{L}_t^\infty(B^\NN_{p,1})})^2\\
&+ (\|{\rm div}\bar{v}\|_{\widetilde{L}_t^\infty(B^\NN_{p,1})}+\|{\rm curl}\bar{v}\|_{\widetilde{L}_t^\infty(B^\NN_{p,1})}+\|{\rm curl}v_L\|_{\widetilde{L}_t^\infty(B^\NN_{p,1})})\|q\|_{\widetilde{L}^1_t(B^{\NN+2}_{p,1})}\\
&\hspace{3cm}+ C[\|q\|_{L_t^\infty(L^\infty)}]\|q\|_{\widetilde{L}_t^1(B^{\NN+2}_{p,1})}(\|{\rm curl}\bar{v}\|_{\widetilde{L}^\infty_t(B^{\NN}_{p,1})}+\|{\rm curl}v_L\|_{\widetilde{L}_t^\infty(B^{\NN}_{p,1})})  \big),
\end{aligned}
\label{estimdiv1}
\end{equation}
Let us now state and prove the following lemma:
\begin{lem}
 \sl{Let $(\bar{q},{\rm div}\bar{v},{\rm curl}\bar{v})$ satisfying (\ref{syst1}) on $[0,T]\times \R^N$. Assume that $(\bar{q},\,{\rm curl}\bar{v},\,q_L,\,{\rm curl}v_L) \in\cC^1([0,T], B^{\NN}_{p,1}\cap B^{\NN+2}_{p,1})$ and ${\rm div}\bar{v}\in\cC^1([0,T],B^{\NN}_{p,1})$, where $q_L$, ${\rm curl}v_L$ satisfy:
\begin{equation}
\begin{cases}
\begin{aligned}
&\p_t q_L-2\mu\D u_l=0,\\
&\p_t{\rm curl v}_L-\mu\D {\rm curl}v_L=0,\\
&(q_L,{\rm curl}v_L)_{t=0}=(q_0,{\rm curl}v_0).
\end{aligned}
\end{cases}
\label{uL}
\end{equation}
There exist a positive constants $\eta$
such that if
  $T$ is small enough so that:
\begin{equation}
\|q_{L}\|_{\widetilde{L}_T^1(B^{\NN+2}_{p,1})}+\|{\rm curl}v_{L}\|_{\widetilde{L}_T^1(B^{\NN+2}_{p,1})}
 \leq \eta^2,
\end{equation}
and $m$ is chosen such that $\|q_0-S_m q_0\|_{\dot{B}_{2,1}^\fd}\leq \eta^2$
then we have, for all $t\in [0,T]$ for $C>0$ large enough:
\begin{equation}
 \begin{cases}
&\|\bar{q}\|_{\widetilde{L}^\infty_T(B^{\NN}_{p,1})}+\|\bar{q}\|_{\widetilde{L}^1_T(B^{\NN+2}_{p,1})}\leq C,\\
&\|{\rm div}v\|_{\widetilde{L}^\infty_T(B^{\NN}_{p,1})}\leq C,\\
&\|{\rm curl}\bar{v}\|_{\widetilde{L}^\infty_T(B^{\NN}_{p,1})}+\||{\rm curl}\bar{v}\|_{\widetilde{L}^1_T(B^{\NN+2}_{p,1})}\leq C
\end{cases}
\end{equation}
}
\label{estimap}
\end{lem}

\textbf{Proof: } 
Assume that $T$ is small enough so that we have for some $\eta\in [0,1]$ (to be precised later) thanks to the classical heat estimates recalled in Proposition \ref{chaleur}:
\begin{equation}
\begin{aligned}
&\|{\rm curl}v_L\|_{\widetilde{L}^1_T (B^{\NN+2}_{p,1})}+\|q_L\|_{\widetilde{L}^1_T (B^{\NN+2}_{p,1})}  \leq \eta,\\
&\|{\rm curl}v_L\|_{\widetilde{L}_T^r (B^{\NN+\frac{2}{r}}_{p,1})}+  \|q_L\|_{\widetilde{L}_T^r (B^{\NN+\frac{2}{r}}_{p,1})}\leq C \eta^{\frac{1}{r}} \|{\rm curl}v_0\|_{B^\NN_{p,1}}^{1-\frac{1}{r}}.
\end{aligned}
\label{reguL2a}
\end{equation}
\begin{equation}
\begin{aligned}
&\|{\rm curl}v_L\|_{\widetilde{L}_t^\infty(B^{\NN}_{p,1})}+\|q_L\|_{\widetilde{L}_t^\infty(B^{\NN}_{p,1})}\leq C (\|{\rm curl}v_0\|_{B^\NN_{p,1}}+\|q_0\|_{B^\NN_{p,1}})=CE_0.
\end{aligned}
\label{reguL1a}
\end{equation}
We are going now to set:
\begin{equation}
\begin{aligned}
&\beta(t)=\|{\rm curl}\bar{v}\|_{\widetilde{L}_t^\infty (B^{\NN}_{p,1})}+\|{\rm curl}\bar{v}\|_{\widetilde{L}_t^1(B^{\NN+2}_{p,1})},\\
&\beta_1(t)=\|q\|_{\widetilde{L}_t^\infty (B^{\NN}_{p,1})}+\|q\|_{\widetilde{L}_t^1(B^{\NN+2}_{p,1})},\\
&\alpha(t)=\|{\rm div}\bar{v}\|_{\widetilde{L}_t^\infty (B^{\NN}_{p,1})}.
\end{aligned}
\end{equation}
and let us define
$$
T^*=\Sup \{t\in[0,T], \quad \beta_1(t) \leq \eta_1,\;\beta(t)\leq \eta_2,\,\;\alpha(t)\leq 2\|{\rm div}v_0\|_{B^{\NN}_{p,1}}=F_0\}.
$$
We recall now that:
$$
\begin{aligned}
\|q\|_{\widetilde{L}^\infty_t(B^\NN_{p,1})}&\leq \beta_1(t)+\|q_L\|_{\widetilde{L}^\infty(B^{\NN}_{p,1})},\\
&\leq  \beta_1(t)+C \|q_0\|_{B^{\NN}_{p,1}}.
\end{aligned}
$$
In the sequel we shall define $\gamma(t)=  \beta_1(t)+C \|q_0\|_{B^{\NN}_{p,1}}$.
Then for all $t\in[0,T^*[$,we have ($\eta\leq 1$) from(\ref{estimq2}),  (\ref{estimdiv1}) and (\ref{estimcurl}):
\begin{equation}
\begin{aligned}
&\beta_1(t) \leq  C\biggl(\sqrt{t}(\beta_1(t)+\|q_L\|_{\widetilde{L}_t^2(B^{\NN+1}_{p,1})})\biggl(\|v_0\|_{L^\infty}+(1+C[\gamma(t)]\gamma(t))\\
&\times \biggl(\sqrt{t}(\beta(t)+ \|{\rm curl} v_L\|_{\widetilde{L}^2(B^{\NN+1}_{p,1})})C[1+\gamma(t)]\\
&+(\beta(t)+\|{\rm curl}v_L\|_{\widetilde{L}_t^{\infty}(B^{\NN}_{p,1})})\big(1+(\sqrt{t}\beta_1(t)+\sqrt{t} \|q_L\|_{\widetilde{L}_1^2(B^{\NN+1}_{p,1})})C[\gamma(t)]\big)\biggl)\\
&+C[\gamma(t)]\sqrt{t}(\beta_1(t)+\|q_L\|_{\widetilde{L}^2(B^{\NN+1}_{p,1})})\biggl)\\
&\hspace{1cm}+
t (\beta(t)+\alpha(t))\gamma(t)  +t(\beta(t)+\alpha(t))
\biggl).
\end{aligned}
\label{estimq2a}
\end{equation}
In a similar way we have:
\begin{equation}
\begin{aligned}
&\beta(t) \leq Ce^{C (t \beta(t) +t \|{\rm curl}v_L\|_{\widetilde{L}_t^\infty(B^{\NN}_{p,1})}+t\alpha(t)+\beta_1(t)+\| q_L\|_{\widetilde{L}_t^1(B^{\NN+2}_{p,1})} +\beta_1(t)+\| q_L\|_{\widetilde{L}_t^2(B^{\NN+1}_{p,1})}    }\biggl(\\
& C[\gamma(t)](\beta_1(t)+\|q_L\|_{\widetilde{L}_t^2(B^{\NN+1}_{p,1})})(\beta(t)+\|{\rm curl}v_L\|_{\widetilde{L}^2_t(B^{\NN+1}_{p,1})})   \\
&+C[\gamma(t)]\big( \|(Id-\dot{S}_m)q_0\|_{B^\NN_{p,1}}\beta(t) +(e^{Ct(\beta(t)+\alpha(t)+\|{\rm curl}v_L\|_{\widetilde{L}^\infty_t(B^{\NN}_{p,1})})}-1)(1+\|q_0\|_{B^{\NN}_{p,1}})        \big)\beta(t)\
\\
&C[\gamma(t)](\beta_1(t)+\|q_L\|_{\widetilde{L}_t^{\infty}(B^\NN_{p,1})})\|{\rm curl}v_L\|_{\widetilde{L}_t^1(B^{\NN+2}_{p,1})}\\
&+(\sqrt{t}\alpha(t)+\sqrt{t}\beta(t)+\sqrt{t}\|{\rm curl}v_L\|_{\widetilde{L}^\infty(B^\NN_{p,1})})^2\\
&+C[\gamma(t)](\beta_1(t)+\|q_L\|_{\widetilde{L}_t^1(B^{\NN+2}_{p,1})})(\beta(t)
+ \|{\rm curl} v_L\|_{\widetilde{L}_t^\infty(B^{\NN}_{p,1})})
\\
&+ (1+C[\gamma(t) ](\beta_1(t)+\|q_L\|_{\widetilde{L}_t^\infty (B^\NN_{p,1})}))(\beta_1(t)+\|q_L\|_{\widetilde{L}_t^2(B^{\NN+1}_{p,1})})(\beta(t)
+ \|{\rm curl} v_L\|_{\widetilde{L}_t^2(B^{\NN+1}_{p,1})})\biggl).
\end{aligned}
\label{estimcurlabc}
\end{equation}
and similarly:
\begin{equation}
\begin{aligned}
&\alpha(t)
\leq Ce^{C (t\alpha(t)+t\beta(t)+t\|{\rm curl}v_L\|_{\widetilde{L}_t^\infty(B^{\NN}_{p,1})}+\beta_1(t)+ \|q_L\|_{\widetilde{L}_t^1(B^{\NN+2}_{p,1})} )     }\big(\|{\rm div}v_0\|_{B^\NN_{p,1}}\\
&+(1+C[\gamma(t)]\gamma(t) )(\beta_1(t)+\|q_L\|_{\widetilde{L}_t^2(B^{\NN+1}_{p,1})})(\beta(t)
+ \|{\rm curl} v_L\|_{\widetilde{L}_t^2(B^{\NN+1}_{p,1})})\\
&+C[\gamma(t)](\beta_1(t)+\|q_L\|_{\widetilde{L}_t^1(B^{\NN+2}_{p,1})})+ (\sqrt{t}\alpha(t)+\sqrt{t}\beta(t)+\sqrt{t}\|{\rm curl}v_L\|_{\widetilde{L}_t^\infty(B^\NN_{p,1})})^2\\
&+(\alpha(t)+\beta(t)+\|{\rm curl}v_L\|_{\widetilde{L}^\infty_t(B^{\NN}_{p,1})})(\beta_1(t)+\|q_L\|_{\widetilde{L}^1_t(B^{\NN+2}_{p,1})})\\
&\hspace{3cm}+ C[\gamma(t)](\beta_1(t)+\|q_L\|_{\widetilde{L}_t^1(B^{\NN+2}_{p,1})})(\beta(t)+\|{\rm curl}v_L\|_{\widetilde{L}^\infty(B^{\NN}_{p,1})})  \big),
\end{aligned}
\label{aestimdiv1b}
\end{equation}
Now using (\ref{reguL1}) and  (\ref{reguL2})  we have then:
\begin{equation}
\begin{aligned}
&\beta_1(t) \leq C\biggl(\sqrt{t}(\beta_1(t)+C\sqrt{\eta E_0})\biggl(\|v_0\|_{L^\infty}+(1+C[\gamma(t)]\gamma(t))\\
&\times \biggl(\sqrt{t}(\beta(t)+C\sqrt{\eta E_0})C[1+\gamma(t)]+(\beta(t)+CE_0)\big(1+(\sqrt{t}\beta_1(t)+\sqrt{t}C\sqrt{\eta E_0})C[\gamma(t)]\big)\biggl)\\
&+C[\gamma(t)]\sqrt{t}(\beta_1(t)+C\sqrt{\eta E_0})\biggl)+
t (\beta(t)+\alpha(t))\gamma(t)  +t(\beta(t)+\alpha(t))
\biggl).
\end{aligned}
\label{estimq2aa}
\end{equation}
It gives in particular since $t<T^{*}$ and $\gamma(t)\leq\eta_1+CE_0$:
\begin{equation}
\begin{aligned}
&\beta_1(t) \leq C\biggl(\sqrt{t}(\eta_1+C\sqrt{\eta E_0})\biggl(\|v_0\|_{L^\infty}+(1+C[\eta_1+CE_0](\eta_1+CE_0))\\
&\times \biggl(\sqrt{t}(\eta_2+C\sqrt{\eta E_0})C[1+\eta_1+CE_0]+(\eta_2+CE_0)\big(1+(\sqrt{t}\eta_1+\sqrt{t}C\sqrt{\eta E_0})C[\eta_1+CE_0]\big)\biggl)\\
&+C[\eta_1+CE_0]\sqrt{t}(\eta_1+C\sqrt{\eta E_0})\biggl)+
t (\eta_2+2F_0)(\eta_1+CE_0)  +t(\eta_2+2F_0)
\biggl).
\end{aligned}
\label{estimq2aaa}
\end{equation}
Next we must choose $t$ small enough in function of $\eta$, $\eta_1$ and $\eta_2$ such that:
\begin{equation}
\begin{aligned}
& C\sqrt{t}\biggl(\|v_0\|_{L^\infty}+(1+C[\eta_1+CE_0](\eta_1+CE_0))\\
&\times \biggl(\sqrt{t}(\eta_2+C\sqrt{\eta E_0})C[1+\eta_1+CE_0]+(\eta_2+CE_0)\big(1+(\sqrt{t}\eta_1+\sqrt{t}C\sqrt{\eta E_0})C[\eta_1+CE_0]\big)\biggl)\\
&+C[\eta_1+CE_0]\sqrt{t}(\eta_1+C\sqrt{\eta E_0})\biggl)\leq\frac{1}{16},\\[3mm]
& C^2\sqrt{\eta E_0}\sqrt{t}\biggl(\|v_0\|_{L^\infty}+(1+C[\eta_1+CE_0](\eta_1+CE_0))\\
&\times \biggl(\sqrt{t}(\eta_2+C\sqrt{\eta E_0})C[1+\eta_1+CE_0]+(\eta_2+CE_0)\big(1+(\sqrt{t}\eta_1+\sqrt{t}C\sqrt{\eta E_0})C[\eta_1+CE_0]\big)\biggl)\\
&+C[\eta_1+CE_0]\sqrt{t}(\eta_1+C\sqrt{\eta E_0})\biggl)\leq\frac{1}{16}\eta_1,\\[3mm]
&t (\eta_2+2F_0)\leq\frac{1}{16},\\[3mm]
&t (\eta_2+2F_0)(1+CE_0)\leq\frac{1}{16}\eta_1.
\end{aligned}
\label{excondi1}
\end{equation}
With the assumption (\ref{excondi1}) we show that:
\begin{equation}
\beta_1(t)\leq\frac{1}{4}\eta_1,
\label{dens1}
\end{equation}
for any $t$ verifying (\ref{excondi1}).
In a similar way we have for $t\leq T^{*}$ using (\ref{estimcurlabc}):
\begin{equation}
\begin{aligned}
&\beta(t) \leq Ce^{C (t \eta_2 +t CE_0+2\eta_1+t F_0+\eta + C\sqrt{E_0\eta} )   }\biggl(\\
& C[  \eta_1+CE_0](\eta_1+C\sqrt{E_0\eta})(\eta_2+C\sqrt{E_0\eta})   \\
&+C[ \eta_1+CE_0 ]\big( \|(Id-\dot{S}_m)q_0\|_{B^\NN_{p,1}}+(e^{Ct(\eta_2+F_0+CE_0)}-1)(1+\|q_0\|_{B^{\NN}_{p,1}})        \big)\eta_2\\
&+C[ \eta_1+CE_0](\eta_1+CE_0)\eta\\
&+(\sqrt{t}F_0+\sqrt{t}\eta_2+\sqrt{t}CE_0)^2+C[\eta_1+CE_0](\eta_1+\eta)(\eta_2
+CE_0)
\\
&+ (1+C[\eta_1+CE_0 ](\eta_1+CE_0))(\eta_1+C\sqrt{E_0\eta})(\eta_2
+C\sqrt{E_0\eta})\biggl).
\end{aligned}
\label{estimcurlaad}
\end{equation}
It gives in particular:
\begin{equation}
\begin{aligned}
&\beta(t) \leq \eta_2 Ce^{C (t \eta_2 +t CE_0+2\eta_1+t F_0+\eta + C\sqrt{E_0\eta} )   }\biggl(\\
& C[  \eta_1+CE_0](\eta_1+C\sqrt{E_0\eta})  \\
&+C[ \eta_1+CE_0 ]\big( \|(Id-\dot{S}_m)q_0\|_{B^\NN_{p,1}}+(e^{Ct(\eta_2+F_0+CE_0)}-1)(1+\|q_0\|_{B^{\NN}_{p,1}})        \big)\\
&+t\eta_2+2\sqrt{t}(\sqrt{t}F_0+\sqrt{t}CE_0)+C[\eta_1+CE_0](\eta_1+\eta)
\\
&+ (1+C[\eta_1+CE_0 ](\eta_1+CE_0))(\eta_1+C\sqrt{E_0\eta})\biggl)\\[3mm]
&+Ce^{C (t \eta_2 +t CE_0+2\eta_1+t F_0+\eta + C\sqrt{E_0\eta} )   }\biggl(\\
& C[  \eta_1+CE_0](\eta_1+C\sqrt{E_0\eta})C\sqrt{E_0\eta} +C[ \eta_1+CE_0](\eta_1+CE_0)\eta\\
&+t(F_0+CE_0)^2+C[\eta_1+CE_0](\eta_1+\eta)CE_0
\\
&+ (1+C[\eta_1+CE_0 ](\eta_1+CE_0))(\eta_1+C\sqrt{E_0\eta})C\sqrt{E_0\eta}\biggl).
\end{aligned}
\label{estimcurlaadde}
\end{equation}
We have now  to choose $t$, $\eta_1$ and $\eta$ small enough and $m$ large enough such that:
\begin{equation}
\begin{aligned}
& Ce^{C (t \eta_2 +t CE_0+2\eta_1+t F_0+\eta + C\sqrt{E_0\eta} )   }\biggl(\\
& C[  \eta_1+CE_0](\eta_1+C\sqrt{E_0\eta})  \\
&+C[ \eta_1+CE_0 ]\big( \|(Id-\dot{S}_m)q_0\|_{B^\NN_{p,1}}+(e^{Ct(\eta_2+F_0+CE_0)}-1)(1+\|q_0\|_{B^{\NN}_{p,1}})        \big)\\
&+t\eta_2+2\sqrt{t}(\sqrt{t}F_0+\sqrt{t}CE_0)+C[\eta_1+CE_0](\eta_1+\eta)
\\
&+ (1+C[\eta_1+CE_0 ](\eta_1+CE_0))(\eta_1+C\sqrt{E_0\eta})\biggl)\leq\frac{1}{4}\\[3mm]
&+Ce^{C (t \eta_2 +t CE_0+2\eta_1+t F_0+\eta + C\sqrt{E_0\eta} )   }\biggl(\\
& C[  \eta_1+CE_0](\eta_1+C\sqrt{E_0\eta})C\sqrt{E_0\eta} +C[ \eta_1+CE_0](\eta_1+CE_0)\eta\\
&+t(F_0+CE_0)^2+C[\eta_1+CE_0](\eta_1+\eta)CE_0
\\
&+ (1+C[\eta_1+CE_0 ](\eta_1+CE_0))(\eta_1+C\sqrt{E_0\eta})C\sqrt{E_0\eta}\biggl)\leq\frac{\eta_2}{4}.
\end{aligned}
\label{extcondi2}
\end{equation}
In a similar way we have using (\ref{aestimdiv1b}) :
\begin{equation}
\begin{aligned}
&\alpha(t)
\leq e^{C (t F_0+t\eta_2+tCE_0+\eta_1+\eta)     }\biggl(\|{\rm div}v_0\|_{B^\NN_{p,1}}\\
&+(1+C[ \eta_1+CE_0] (\eta_1+CE_0) )(\eta_1+C\sqrt{E_0 \eta})(\eta_2
+ C\sqrt{E_0 \eta})\\
&+C[\eta_1+CE_0](\eta_1+C\sqrt{E_0 \eta})+ (\sqrt{t}F_0+\sqrt{t}\eta_2+\sqrt{t}CE_0)^2\\
&+(F_0+\eta_2+CE_0)(\eta_1+\eta)+ C[\eta_1+CE_0] (\eta_1+\eta)(\eta_2+CE_0)  \biggl),
\end{aligned}
\label{estimdiv1aae}
\end{equation}
It suffices now to choose $\eta$, $\eta_1$, $\eta_2$ and the time $t$ small enough such that:
\begin{equation}
\begin{aligned}
&e^{C (t F_0+t\eta_2+tCE_0+\eta_1+\eta)     }\leq\frac{4}{3},\\[2mm]
 & e^{C (t F_0+t\eta_2+tCE_0+\eta_1+\eta)     } \biggl((1+C[ \eta_1+CE_0] (\eta_1+CE_0) )(\eta_1+C\sqrt{E_0 \eta})(\eta_2
+ C\sqrt{E_0 \eta})\\
&+C[\eta_1+CE_0](\eta_1+C\sqrt{E_0 \eta})+ (\sqrt{t}F_0+\sqrt{t}\eta_2+\sqrt{t}CE_0)^2\\
&+(F_0+\eta_2+CE_0)(\eta_1+\eta)+ C[\eta_1+CE_0] (\eta_1+\eta)(\eta_2+CE_0)  \biggl)\leq \frac{1}{3}\|{\rm div}v_0\|_{B^{\NN}_{p,1}}.
\end{aligned}
\label{extcondi3}
\end{equation}
Choosing $T$, $\eta$, $\eta_1$ and $\eta_2$ as in conditions (\ref{reguL2a}), (\ref{excondi1}), (\ref{extcondi2}) and (\ref{extcondi3}) allows to ensure that  $T^*\geq T>0$. An it concludes the proof of the lemma.
\subsection{Existence}
We use a standard scheme for proving the existence of the solutions:
\begin{enumerate}
\item We smooth out the data and get a sequence of smooth solutions $(q^{n},u^{n})_{n\in\mathbb{N}}$ of an approximated system of (\ref{0.6}), on a bounded interval $[0,T^{n}]$ which may depend on $n$. 
\item We exhibit a positive lower bound $T$ for $T^{n}$, and prove uniform estimates on $(\bar{q}^{n},{\rm div}\bar{v}^{n},{\rm curl}\bar{v}^n)$ (we refer to the next subsection for the definitions ) in the space
\begin{equation}
 E_{T}=\big(\widetilde{C}_{T}(B^{\NN}_{p,1})\cap \widetilde{L}^1(B^{\NN+2}_{p,1})\big)\times \widetilde{C}_{T}(B^{\NN}_{p,1})\times\big( \widetilde{C}_{T}(B^{\NN}_{p,1})\cap \widetilde{L}^1(B^{\NN+2}_{p,1})\big).
\label{ET}
\end{equation}
\item We use compactness to prove that the sequence $(q^{n},{\rm div}v^{n}, {\rm curl}v^{n}  )$ converges, up to extraction, to a solution of (\ref{0.1}).
\end{enumerate}
\subsubsection*{Step 1: Approximation}
It suffices to regularize the initial data by choosing $(q_0^n, u_0^n)$ belonging to $B^{\NN}_{p,1}\times B^{\NN-1}_{p,1}$ with $v_0^n$ in $L^\infty\cap B^{\NN+1}_{p,1}$. We know that it exists $T_n$ such that we have a strong solution in $E_T$ (see \cite{M3AS}). Using the lemma \ref{estimap} we prove that $T_n\geq T$ since we control the Lipschitz norm $\n u^n$  in $L^1(L^\infty)$ (see \cite{M3AS}).
\subsubsection*{Step 2: compactness and convergence}

This part is also classical and we refer for example to \cite{DW} (chapter $10$) for details: using the previous result and the Ascoli theorem, we can extract a subsequence that weakly converges towards some couple $(q,v)$, which is proved to be a solution of the original system and to satisfy the energy estimates. This concludes the existence part of the theorem.

\section{Uniqueness}
\label{section4}

We are going now to prove the uniqueness in the following space:
$$
F_T=\big(\widetilde{L}^\infty_T(B^{\NN}_{p,1})\cap \widetilde{L}^1_T(B^{\NN+2}_{p,1})\big)\times \widetilde{L}^\infty_T(B^{\NN}_{p,1})\times \big(\widetilde{L}^\infty_T(B^{\NN}_{p,1})\cap \widetilde{L}^1_T(B^{\NN+2}_{p,1})\big).
$$

\begin{theorem}
 \sl{Let $N\geq 2$ and assume that $(q_i, v_i)$ ($i\in\{1,2\}$) are two solutions of  (\ref{0.6}) with the same initial data on the same interval $[0, T]$ and both belonging to the space $F_T$ then $(q_1,u_1)\equiv(q_2,u_2)$ on $[0,T]$.
}
\end{theorem}

\textbf{Proof: } for $i\in\{1,2\}$, let us recall that $(\bar{q}_i, \bar{v}_i)$ satisfy the system:

\begin{equation}
\begin{cases}
\begin{aligned}
&\p_t \bar{q}_i-2\mu'(1+q_1)\D q_i+v_i\cdot\n \bar{q}_i
=F_i,\\
&\p_t {\rm div}\bar{v}_i+ u_i\cdot\n {\rm div}\bar{v}_i=G_i,\\
&\p_t {\rm curl}\bar{v}_i+u_i\cdot\n{\rm curl}v_i-\big(\mu+\dot{S}_m\big(\frac{\mu(1+q_i)}{1+q_i}-\mu(1)\big)\big)\D{\rm curl}\bar{v}_i=H_i.
\end{aligned}
\end{cases}
\label{syst1a}
\end{equation}
with:
$$
\begin{aligned}
&F_i=-(q_i+1){\rm div}v_i-v_i\cdot\n q_L+2\mu''(1+q_i)|\n q_i|^2,\\
&G_i=-\n v_i: ^{t}\n u_i+\frac{1}{2}\frac{\lambda(1+q_i)}{(1+q_i)^2}\n q_i\cdot{\rm div}{\rm curl}v_i+R(\va_i,v_i)+\frac{1}{2}\n\n\va(1+q_i):{\rm curl}v_i -\D F(1+q_i),\\
&H_i=\frac{1}{2}\n\va(1+q_i)\cdot\n{\rm curl}v_i+(Id-\dot{S}_m)\big(\frac{\mu(1+q_i)}{1+q_i}-\mu\big)\D{\rm curl}\bar{v}_i+(\frac{\mu(1+q_i)}{1+q_i}-\mu)\D{\rm curl}v_L-R^i_1,
\end{aligned}
$$
and
$$
\begin{aligned}
&(R^i_1)_{i'j}=\sum_k (\p_{i'} u^i_k\p_k v^i_j-\p_j u^i_k\p_k v^i_{i'})-\frac{1}{2}\frac{\lambda(\rho^i)}{(\rho^i)^2}\big(\sum_k (\p_{i'}\rho^i\p_k({\rm curl}v^i)_{kj}-\p_j\rho^i\p_k({\rm curl}v^i)_{ki'})\\
&\hspace{7cm}-\sum_k \big(\p_{i'k}\va(\rho^i)({\rm curl}v^i)_{k j}-\p_{k j}\va(\rho^i)({\rm curl}v^i)_{i'k}\big).
\end{aligned}
$$
Let us recall that $(q_L,{\rm div}v_L,{\rm curl}v_L)$ is defined as in the previous section. If we denote by $\dq=q_1-q_2$ and $\delta v=v_1-v_2$, then $(\dq,\delta v)$ satisfy the following system:
\begin{equation}
\begin{cases}
\begin{aligned}
&\p_t \delta q-2\mu'(1+q_1)\D \delta q+v_1\cdot\n \delta q
=\delta F,\\
&\p_t {\rm div}\delta v+ u_1\cdot\n {\rm div}\delta v=\delta G,\\
&\p_t {\rm curl}\delta v+u_1\cdot\n{\rm curl}\delta v-\big(\mu+\dot{S}_m\big(\frac{\mu(1+q_1)}{1+q_1}-\mu(1)\big)\big)\D{\rm curl}\delta v=\delta H.
\end{aligned}
\end{cases}
\label{syst1b}
\end{equation}
with:
$$
\begin{aligned}
&\delta F=\delta F_1+\delta F_2+\delta F_3,\\
&\delta G=\delta G_1+\delta G_2+\delta G_3+\delta G_4,\\
&\delta H=\delta H_1+\delta H_2+\delta H_3+\delta H_4+\delta H_5.
\end{aligned}
$$
and:
$$
\begin{cases}
\de F_1=-2\D q_2(\mu'(1+q_2)-\mu'(1+q_1))-\delta v\cdot\n\bar{q}_2   ,\\
\de F_2=-(1+q_1)\div \delta v-\delta q{\rm div}v_2-\delta v\cdot\n q_L,\\
\de F_3= 2(\mu''(1+q_1)-\mu''(1+q_2) )|\n q_1|^2+2\mu''(1+q_2)\n \delta q\cdot\n(q_1+q_2),\\
\de G_1=-\delta u\cdot\n{\rm div}\bar{v}_2-\n\delta v:^t\n u_1-\n v_2:^{t}\n\delta u,\\
\de G_2= \frac{1}{2}\frac{\lambda(1+q_1)}{(1+q_1)^2}\big( \n \delta q\cdot{\rm div}{\rm curl}v_1+\n q_2\cdot{\rm div}{\rm curl}\delta v\big)\\
\hspace{3cm}+ \frac{1}{2}(\frac{\lambda(1+q_1)}{(1+q_1)^2}- \frac{\lambda(1+q_2)}{(1+q_2)^2}) \n q_2\cdot{\rm div}{\rm curl}v_2   ,\\
\de G_3=-\frac{1}{2}\p_i\va(\rho_1)\p_j({\rm curl}\delta v_{ij})-\frac{1}{2}\p_i(\va(\rho_1)-\va(\rho_2))\p_j (({\rm curl}v_2)_{ij}) \\
\hspace{6cm}-\D \big[ F(1+q_1)-F(1+q_2)],\\
\delta G_4=\frac{1}{2}\n\n(\va(1+q_1)-\va(1+q_2)):{\rm curl}v_1+\n\n\va(1+q_2):{\rm curl}\delta v,\\ 
\de H_1=-\delta u\cdot\n{\rm curl}\bar{v}_2-S_m\big(\frac{\mu(1+q_1)}{1+q_1}-\frac{\mu(1+q_1)}{1+q_1}\big)\D{\rm curl}\bar{v}_2 ,\\
\delta H_2=\frac{1}{2}\big(\n\va(1+q_1)\cdot\n{\rm curl}\delta v+\n(\va(1+q_1)-\va(1+q_2))\cdot\n{\rm curl}v_2\big),\\
\delta H_3=(Id-\dot{S}_m)\big(\frac{\mu(1+q_1)}{1+q_1}-\mu\big)\D{\rm curl}\delta v+
(Id-\dot{S}_m)\big(\frac{\mu(1+q_1)}{1+q_1}-(\frac{\mu(1+q_2)}{1+q_2})\D{\rm curl}\bar{v}_2\\
\delta H_4=(\frac{\mu(1+q_1)}{1+q_1}-\frac{\mu(1+q_2)}{1+q_2})\D{\rm curl}v_L\\
\delta H_5=R_1^1-R_1^2.
\end{cases}
$$
We have in particular after tedious calculus:
$$
\begin{aligned}
&(R_1^1-R_1^2)_{ij}=\sum_k(\p_iu_k^1\p_k\delta v_j-\p_j u_k^1\p_k\delta v_i)+(\p_i\delta u_k\p_k v_j^2-\p_j\delta u_k\p_k v_i^2)\\
&-\frac{1}{2}\frac{\lambda(\rho_1)}{\rho_1^2}\sum_k\big(\p_i\rho_1\p_k({\rm curl}\delta v)_{kj}+\p_i\delta\rho\p_k({\rm curl}v_2)_{kj}\big)\\
&-\frac{1}{2}(\frac{\lambda(\rho_1)}{\rho_1^2}- \frac{\lambda(\rho_2)}{\rho_2^2})\sum_k(\p_i\rho_2\p_k({\rm curl}v_2)_{kj}-\p_j\rho_2\p_k(({\rm curl}v_2)_{ki})\\[2mm]
&+\frac{1}{2}\frac{\lambda(\rho_1)}{\rho_1^2}\sum_k\big(\p_{ik}\va(\rho_1){\rm curl}\delta v_{kj}+\p_{kj}(\va(\rho_1)-\va(\rho_2)) ({\rm curl}v_2)_{ki}\big)\\
&-\frac{1}{2}(\frac{\lambda(\rho_1)}{\rho_1^2}- \frac{\lambda(\rho_2)}{\rho_2^2})\sum_k(\p_{ik}\va(\rho_2)({\rm curl}v_2)_{kj}-\p_{kj}\va(\rho_2)({\rm curl}v_2)_{ki})\\
\end{aligned}
$$
We wish to prove (as for $(NSC)$) the uniqueness in the following space:
$$
E_T=\big(\widetilde{L}^\infty_T(B^{\NN-1}_{p,1})\cap \widetilde{L}^1_T(B^{\NN+1}_{p,1})\big)\times \widetilde{L}^\infty_T(B^{\NN-1}_{p,1})\times \big(\widetilde{L}^\infty_T(B^{\NN-1}_{p,1})\cap \widetilde{L}^1_T(B^{\NN+1}_{p,1})\big).
$$
Due to the term $\delta u\cdot\n{\rm div}\bar{v}_2$ in the right-hand side of the second equation in system (\ref{syst1b}), we loose one
derivative when estimating ${\rm div}\delta v$: one only gets bounds in $\widetilde{L}^\infty_T(B^{\NN-1}_{p,1})$. Now, the right hand
side of the first and the third equations contains a term of type $\delta v \cdot\n \bar{q}_2$ and $\delta u\cdot\n{\rm curl}\bar{v}_2$ so that the loss of one derivative for ${\rm div}\delta v$ entails a loss of one derivative for $\delta q$ and ${\rm curl}\delta v$. Therefore, getting bounds in $E_T$
for $(\delta q,{\rm div}\delta v,{\rm curl}\delta v)$ is the best that one can hope.\\

Let us begin with $\dq$. As $q_1$ and $q_2$ have the same initial data, using the proposition \ref{lemmeuti} leads to:
\begin{equation}
\begin{aligned}
&\|\delta q\|_{\widetilde{L}^\infty_t(B^{\NN-1}_{p,1})}+ \|\delta q\|_{\widetilde{L}^1_t(B^{\NN+1}_{p,1})}\leq e^{C V(t)}\int^t_0\|\delta F(\tau)\|_{B^{\NN-1}_{p,1}}d\tau.
\end{aligned}
\end{equation}
with $V(t)=\int^t_0(\|v_1(\tau)\|_{B^{\NN+1}_{p,1}}+\|\mu'(1+q_1)-\mu'(1)\|^2_{B^{\NN+1}_{p,1}})d\tau$. Using propositions \ref{produit} and \ref{estimcompo}, we have for $p<N$:
$$
\begin{aligned}
&\|\D q_2(\mu'(1+q_2)-\mu'(1+q_1))\|_{B^{\NN-1}_{p,1}}\lesssim \|q_2\|_{B^{\NN+2}_{p,1}}\|\delta q\|_{B^{\NN-1}_{p,1}} C[\|q_1\|_{L^\infty},\|q_2\|_{L^\infty}],\\
&\|\delta v\cdot\n\bar{q}_2\|_{B^{\NN-1}_{p,1}}\lesssim \|\delta v\|_{B^{\NN}_{p,1}}\|\bar{q}_2\|_{B^{\NN}_{p,1}}  ,\\
&\|(1+q_1)\div \delta v\|_{B^{\NN-1}_{p,1}}\lesssim (1+\|q_1\|_{B^{\NN}_{p,1}})\|{\rm div} \delta v\|_{B^{\NN-1}_{p,1}},\\
&\| \delta q{\rm div}v_2\|_{B^{\NN-1}_{p,1}}\lesssim \|\delta q\|_{B^{\NN-1}_{p,1}}\|{\rm div}v_2\|_{B^{\NN}_{p,1}},\\
&\|\delta v\cdot\n q_L\|_{B^{\NN-1}_{p,1}}\lesssim  \|\delta v\|_{B^{\NN}_{p,1}}\| q_L\|_{B^{\NN}_{p,1}},\\
\end{aligned}
$$
$$
\begin{aligned}
&\|(\mu''(1+q_1)-\mu''(1+q_2) )|\n q_1|^2\|_{B^{\NN-1}_{p,1}}\lesssim\|\delta q\|_{B^{\NN-1}_{p,1}}\|q_1\|^2_{B^{\NN+1}_{p,1}}C[\|q_1\|_{L^\infty},\|q_2\|_{L^\infty}],\\
&\|\mu''(1+q_2)\n \delta q\cdot\n(q_1+q_2)\|_{B^{\NN-1}_{p,1}}\lesssim \|\delta q\|_{B^{\NN}_{p,1}}(\|q_1\|_{B^{\NN+1}_{p,1}}+\|q_2\|_{B^{\NN+1}_{p,1}})\\
&\hspace{8cm}\times(1+C[\|q_2\|_{L^\infty}]\|q_2\|_{B^\NN_{p,1}}).
\end{aligned}
$$
Collecting all the previous inequalities, we have:
\begin{equation}
\begin{aligned}
&\|\delta q\|_{\widetilde{L}^\infty_t(B^{\NN-1}_{p,1})}+ \|\delta q\|_{\widetilde{L}^1_t(B^{\NN+1}_{p,1})}\lesssim e^{C V(t)}\int^t_0\big( \|\delta q(\tau)\|_{B^{\NN-1}_{p,1}} (\|q_2\|_{B^{\NN+2}_{p,1}}+\|{\rm div}v_2\|_{B^{\NN}_{p,1}}\\
&+\|q_1\|^2_{B^{\NN+1}_{p,1}})+(\|{\rm div}\delta v\|_{B^{\NN-1}_{p,1}}+ \|{\rm curl}\delta v\|_{B^{\NN-1}_{p,1}})    (1+\|\bar{q}_1\|_{B^{\NN}_{p,1}}+\|\bar{q}_2\|_{B^{\NN}_{p,1}}+2\|q_L\|_{B^{\NN}_{p,1}}) \big)  d\tau\\
&+ e^{C V(t)}\int^t_0\big \|\delta q\|_{B^{\NN}_{p,1}}(\|q_1\|_{B^{\NN+1}_{p,1}}+\|q_2\|_{B^{\NN+1}_{p,1}})(1+C(\|q_2\|_{L^\infty}\|q_2\|_{B^\NN_{p,1}})d\tau.
\end{aligned}
\label{transportuni}
\end{equation}
Let us deal now with the second equation of the system (\ref{syst1b}) and using the proposition \ref{transport1} we have
\begin{equation}
\begin{aligned}
&\|{\rm div}\delta v\|_{\widetilde{L}^\infty_t(B^{\NN-1}_{p,1})}\leq e^{C V_1(t)}\int^t_0\|\delta G(\tau)\|_{B^{\NN-1}_{p,1}}d\tau,
\end{aligned}
\end{equation}
with $V_1(t)=\int^t_0\|u_1(\tau)\|_{B^{\NN+1}_{p,1}}d\tau$. Next using proposition \ref{produit}and \ref{estimcompo} we have:
$$
\begin{aligned}
&\|\delta u\cdot\n{\rm div}\bar{v}_2\|_{B^{\NN-1}_{p,1}}\lesssim \|\delta u\|_{B^\NN_{p,1}}\|{\rm div}\bar{v}_2\|_{B^{\NN}_{p,1}},\\
&\|\n\delta v:^t\n u_1\|_{B^{\NN-1}_{p,1}}\lesssim\|\delta v\|_{B^{\NN}_{p,1}} \|u_1\|_{B^{\NN+1}_{p,1}},\\
&\|\n v_2:^{t}\n\delta u\|_{B^{\NN-1}_{p,1}} \lesssim\|\delta u\|_{B^{\NN}_{p,1}} \|v_2\|_{B^{\NN+1}_{p,1}},\\
&\|\frac{\lambda(1+q_1)}{(1+q_1)^2} \n \delta q\cdot{\rm div}{\rm curl}v_1\|_{B^{\NN-1}_{p,1}} \lesssim \| \delta q\|_{B^{\NN}_{p,1}} \|{\rm div}{\rm curl}v_1\|_{B^{\NN}_{p,1}}\|q_1\|_{B^\NN_{p,1}}C[\|q_1\|_{L^\infty}]  ,\\ 
&\| \n q_2\cdot{\rm div}{\rm curl}\delta v\|_{B^{\NN-1}_{p,1}}\lesssim \| {\rm curl} \delta v\|_{B^{\NN}_{p,1}}\|q_2\|_{B^{\NN+1}_{p,1}},\\
&\| (\frac{\lambda(1+q_1)}{(1+q_1)^2}- \frac{\lambda(1+q_2)}{(1+q_2)^2}) \n q_2\cdot{\rm div}{\rm curl}v_2 \|_{B^{\NN-1}_{p,1}}\lesssim \|\delta q\|_{B^{\NN-1}_{p,1}}      \|q_2\|_{B^{\NN+1}_{p,1}}  \\
&\hspace{7cm}\times \|{\rm curl}v_2 \|_{B^{\NN+1}_{p,1}} C[\|q_1\|_{L^\infty}, \|q_2\|_{L^\infty}]   ,\\
&\| \p_i\va(\rho_1)\p_j({\rm curl}\delta v_{ij})\|_{B^{\NN-1}_{p,1}}\lesssim \|q_1\|_{B^{\NN+1}_{p,1}} C[\|q_1\|_{L^\infty}] \|{\rm curl}\delta v\|_{B^{\NN}_{p,1}},\\
&\|\p_i(\va(\rho_1)-\va(\rho_2))\p_j(({\rm curl}v_2)_{ij}) \|_{B^{\NN-1}_{p,1}}\lesssim \|\delta q\|_{B^{\NN}_{p,1}} C[\|q_1\|_{L^\infty}, \|q_2\|_{L^\infty}]  \| {\rm curl}v_2\|_{B^{\NN+1}_{p,1}},\\
&\|\D \big[ F(1+q_1)-F(1+q_2)]\|_{B^{\NN-1}_{p,1}}\lesssim \|\delta q\|_{B^{\NN+1}_{p,1}}C[\|q_1\|_{L^\infty}, \|q_2\|_{L^\infty}]    ,\\
&\|\n\n(\va(1+q_1)-\va(1+q_2)):{\rm curl}v_1\|_{B^{\NN-1}_{p,1}}\lesssim \|\delta q\|_{B^{\NN+1}_{p,1}}C[\|q_1\|_{L^\infty}, \|q_2\|_{L^\infty}] \|{\rm curl}v_1\|_{B^{\NN}_{p,1}},\\
&\| \n\n\va(1+q_2):{\rm curl}\delta v\|_{B^{\NN-1}_{p,1}}\lesssim C[\|q_2\|_{L^\infty}] \|q_2\|_{B^{\NN+2}_{p,1}}\|{\rm curl}\delta v\|_{B^{\NN-1}_{p,1}}.
\end{aligned}
$$
Collecting all the inequalities, we have:
\begin{equation}
\begin{aligned}
&\|{\rm div}\delta v\|_{\widetilde{L}^\infty_t(B^{\NN-1}_{p,1})}\lesssim e^{C V_1(t)}\int^t_0\|\delta v\|_{B^{\NN}_{p,1}} \|u_1(\tau)\|_{B^{\NN+1}_{p,1}}  d\tau,\\
&+ e^{C V_1(t)}\int^t_0  \|\delta u\|_{B^\NN_{p,1}}(\|{\rm div}\bar{v}_2\|_{B^{\NN}_{p,1}}+ \|v_2\|_{B^{\NN+1}_{p,1}})\\
&+e^{C V_1(t)}\int^t_0  \|\delta q\|_{B^{\NN-1}_{p,1}} \|q_2\|_{B^{\NN+1}_{p,1}}   \|{\rm curl}v_2 \|_{B^{\NN+1}_{p,1}} d\tau\\
&+e^{C V_1(t)}\int^t_0 \| \delta q\|_{B^{\NN}_{p,1}}( \|{\rm curl}v_1\|_{B^{\NN+1}_{p,1}}\|q_1\|_{B^\NN_{p,1}}C[\|q_1\|_{L^\infty}] + \| {\rm curl}v_2\|_{B^{\NN+1}_{p,1}})d\tau\\
&+e^{C V_1(t)}\int^t_0  \|\delta q\|_{B^{\NN+1}_{p,1}}\big( C[\|q_1\|_{L^\infty}, \|q_2\|_{L^\infty}] +C[\|q_1\|_{L^\infty}, \|q_2\|_{L^\infty}] \|{\rm curl}v_1\|_{B^{\NN}_{p,1}})d\tau\\
&+e^{C V_1(t)}\int^t_0 \| {\rm curl} \delta v\|_{B^{\NN}_{p,1}}(\|q_2\|_{B^{\NN+1}_{p,1}}+\|q_1\|_{B^{\NN+1}_{p,1}} C[\|q\|_{L^\infty}])d\tau\\
&+e^{C V_1(t)}\int^t_0 \| {\rm curl} \delta v\|_{B^{\NN-1}_{p,1}} C[\|q_2\|_{L^\infty}] \|q_2\|_{B^{\NN+2}_{p,1}}d\tau
\end{aligned}
\label{unidiv}
\end{equation}
We proceed similarly for the third equation of system (\ref{syst1b}) and we apply the proposition \ref{lemmeuti}:
\begin{equation}
\begin{aligned}
&\|{\rm curl}\delta v\|_{\widetilde{L}^\infty_t(B^{\NN-1}_{p,1})}\leq e^{C V_2(t)}\int^t_0\|\delta H\|_{B^{\NN-1}_{p,1}}d\tau.
\end{aligned}
\label{unicurl}
\end{equation}
with $V_2(t)=\int^t_0(\|u_1(\tau)\|_{B^{\NN+1}_{p,1}}+\|\frac{\mu(1+q_1)}{1+q_1}-\mu(1)\|^2_{B^{\NN+1}_{p,1}})d\tau$. We have using proposition \ref{produit} and (\ref{estimcompo}):
$$
\begin{aligned}
&\|\delta u\cdot\n{\rm curl}\bar{v}_2\|_{B^{\NN-1}_{p,1}} \lesssim \|\delta u\|_{B^{\NN}_{p,1}}\|{\rm curl}\bar{v}_2\|_{B^{\NN}_{p,1}},  \\
&\|S_m\big(\frac{\mu(1+q_1)}{1+q_1}-\frac{\mu(1+q_1)}{1+q_1}\big)\D{\rm curl}\bar{v}_2\|_{B^{\NN-1}_{p,1}} \lesssim \|\delta q\|_{B^{\NN-1}_{p,1}} \|{\rm curl}\bar{v}_2\|_{B^{\NN+2}_{p,1}} C[\|q_1\|_{L^\infty},\|q_2\|_{L^\infty}],\\
&\|\n\va(1+q_1)\cdot\n{\rm curl}\delta v\|_{B^{\NN-1}_{p,1}} \lesssim \|q_1\|_{B^{\NN+1}_{p,1}}C[\|q_1\|_{L^\infty}]  \|{\rm curl}\delta v\|_{B^{\NN}_{p,1}} ,\\
&\|\n(\va(1+q_1)-\va(1+q_2))\cdot\n{\rm curl}v_2\|_{B^{\NN-1}_{p,1}} \lesssim\|\delta q\|_{B^{\NN}_{p,1}}C[\|q_1\|_{L^\infty},\|q_2\|_{L^\infty}] \| {\rm curl}v_2\|_{B^{\NN+1}_{p,1}},   \\
&\|(Id-\dot{S}_m)\big(\frac{\mu(1+q_1)}{1+q_1}-\mu\big)\D{\rm curl}\delta v\|_{B^{\NN-1}_{p,1}} \lesssim C[\|q_1\|_{L^\infty}]   \|( Id-\dot{S}_m)q_1)\|_{B^\NN_{p,1}} \|{\rm curl}\delta v\|_{B^{\NN+1}_{p,1}}  \\
\end{aligned}
$$
$$
\begin{aligned}
&\|(Id-\dot{S}_m)\big(\frac{\mu(1+q_1)}{1+q_1}-(\frac{\mu(1+q_2)}{1+q_2})\D{\rm curl}\bar{v}_2\|_{B^{\NN-1}_{p,1}} \lesssim\|\delta q\| _{B^{\NN-1}_{p,1}}C[\|q_1\|_{L^\infty}, \|q_2\|_{L^\infty}]  \|{\rm curl}\bar{v}_2\|_{B^{\NN+2}_{p,1}}  \\
&\|(\frac{\mu(1+q_1)}{1+q_1}-\frac{\mu(1+q_2)}{1+q_2})\D{\rm curl}v_L\|_{B^{\NN-1}_{p,1}} \lesssim
\|\delta q\| _{B^{\NN-1}_{p,1}}C[\|q_1\|_{L^\infty}, \|q_2\|_{L^\infty}]  \|{\rm curl}v_L\|_{B^{\NN+2}_{p,1}} \\[3mm]
&\|R_1^1-R_1^2\|_{B^{\NN-1}_{p,1}} \lesssim\|\delta v\|_{B^{\NN}_{p,1}}\|u_1\|_{B^{\NN+1}_{p,1}}+ 
\|\delta u\|_{B^{\NN}_{p,1}}\|u_2\|_{B^{\NN+1}_{p,1}}\\
&+(1+\|q_1\|_{B^\NN_{p,1}}C[\|q_1\|_{L^\infty}])
\big(\|{\rm curl}\delta v\|_{B^{\NN}_{p,1}}\|q_1\|_{B^{\NN+1}_{p,1}}+\|\delta q\|_{B^{\NN}_{p,1}}\|{\rm curl}v_2\| _{B^{\NN+1}_{p,1}}\big)\\
&+\|\delta q\|_{B^{\NN-1}_{p,1}}\|{\rm curl}v_2\|_{B^{\NN+1}_{p,1}}\|q_2\|_{B^{\NN+1}_{p,1}}+(1+\|q_1\|_{B^{\NN}_{p,1}}C[\|q_1\|_{L^\infty}])\\
&\times\big(\|{\rm curl}\delta v\|_{B^{\NN-1}_{p,1}}\|q_1\| _{B^{\NN+2}_{p,1}}+\|\delta q\|_{B^{\NN+1}_{p,1}}\|{\rm curl}v_2\|_{B^{\NN}_{p,1}}C[\|q_1\|_{L^\infty}, 
\|q_2\|_{L^\infty}]\big)\\
&+\|\delta q\|_{B^{\NN-1}_{p,1}}\|q_2\|_{B^{\NN+2}_{p,1}} C[\|q_2\|_{L^\infty}]\|{\rm curl}v_2\|_{B^{\NN}_{p,1}}.
\end{aligned}
$$
Combining all the previous inequalities we have:
\begin{equation}
\begin{aligned}
&\|{\rm curl}\delta v\|_{\widetilde{L}^\infty_t(B^{\NN-1}_{p,1})}+\|{\rm curl}\delta v\|_{\widetilde{L}^1_t(B^{\NN+1}_{p,1})}\lesssim e^{C V_2(t)}\int^t_0 \|\delta u\|_{B^{\NN}_{p,1}}\|{\rm curl}\bar{v}_2\|_{B^{\NN}_{p,1}}d\tau\\
&+e^{C V_2(t)}\int^t_0  \|\delta q\|_{B^{\NN-1}_{p,1}} \big( \|{\rm curl}\bar{v}_2\|_{B^{\NN+2}_{p,1}}+C[\|q_1\|_{L^\infty}, \|q_2\|_{L^\infty}]  \|{\rm curl}\bar{v}_2\|_{B^{\NN+2}_{p,1}}\\
&+C[\|q_1\|_{L^\infty}, \|q_2\|_{L^\infty}]  \|{\rm curl}v_L\|_{B^{\NN+2}_{p,1}}+\|{\rm curl}v_2\|_{B^{\NN+1}_{p,1}}\|q_2\|_{B^{\NN+1}_{p,1}}\\
&+\|q_2\|_{B^{\NN+2}_{p,1}} C[\|q_2\|_{L^\infty}]\|{\rm curl}v_2\|_{B^{\NN}_{p,1}}\big)d\tau\\
&+e^{C V_2(t)}\int^t_0  \|{\rm curl}\delta v\|_{B^{\NN}_{p,1}}  \|q_1\|_{B^{\NN+1}_{p,1}}C[\|q_1\|_{L^\infty}]d\tau\\
&+e^{C V_2(t)}\int^t_0 \|\delta q\|_{B^{\NN}_{p,1}}\big( C[\|q_1\|_{L^\infty},\|q_2\|_{L^\infty}] \| {\rm curl}v_2\|_{B^{\NN+1}_{p,1}}\\
&+(1+\|q_1\|_{B^\NN_{p,1}}C[\|q_1\|_{L^\infty}])\|{\rm curl}v_2\| _{B^{\NN+1}_{p,1}}\big)d\tau\\
&+e^{C V_2(t)}\int^t_0 \|{\rm curl}\delta v\|_{B^{\NN+1}_{p,1}} C[\|q_1\|_{L^\infty}]   \|( Id-\dot{S}_m)q_1)\|_{B^\NN_{p,1}} d\tau\\
&+e^{C V_2(t)}\int^t_0 (\|\delta v\|_{B^{\NN}_{p,1}}\|u_1\|_{B^{\NN+1}_{p,1}}+ 
\|\delta u\|_{B^{\NN}_{p,1}}\|u_2\|_{B^{\NN+1}_{p,1}})d\tau\\
&+e^{C V_2(t)}\int^t_0 \|{\rm curl}\delta v\|_{B^{\NN-1}_{p,1}} (1+\|q_1\|_{B^{\NN}_{p,1}}C[\|q_1\|_{L^\infty}])\|q_1\| _{B^{\NN+2}_{p,1}}\\
&+e^{C V_2(t)}\int^t_0\|\delta q\|_{B^{\NN+1}_{p,1}} (1+\|q_1\|_{B^{\NN}_{p,1}}C[\|q_1\|_{L^\infty}])\|{\rm curl}v_2\|_{B^{\NN}_{p,1}}C[\|q_1\|_{L^\infty}, 
\|q_2\|_{L^\infty}]d\tau
\end{aligned}
\label{unicurl1}
\end{equation}
From (\ref{transportuni}) we have: 
\begin{equation}
\begin{aligned}
&\|\delta q\|_{\widetilde{L}^\infty_t(B^{\NN-1}_{p,1})}+ \|\delta q\|_{\widetilde{L}^1_t(B^{\NN+1}_{p,1})}\leq C e^{C V(t)}\|\delta q\|_{\widetilde{L}^\infty_t(B^{\NN-1}_{p,1})} \int^t_0\big( 
 (\|q_2\|_{B^{\NN+2}_{p,1}}+\|{\rm div}v_2\|_{B^{\NN}_{p,1}}\\
&+\|q_1\|^2_{B^{\NN+1}_{p,1}})+(\|{\rm div}\delta v\|_{B^{\NN-1}_{p,1}}+ \|{\rm curl}\delta v\|_{B^{\NN-1}_{p,1}})    (1+\|\bar{q}_1\|_{B^{\NN}_{p,1}}+\|\bar{q}_2\|_{B^{\NN}_{p,1}}+2\|q_L\|_{B^{\NN}_{p,1}}  )  d\tau\\
&+Ce^{C V(t)} \|\delta q\|_{\widetilde{L}^2(B^{\NN}_{p,1})}(\|q_1\|_{L^2(B^{\NN+1}_{p,1})}+\|q_2\|_{L^2(B^{\NN+1}_{p,1})})(1+C(\|q_2\|_{L^\infty})\|q_2\|_{\widetilde{L}^\infty(B^\NN_{p,1})}).
\end{aligned}
\label{transportuni1}
\end{equation}
It implies by Young inequality and $C$ large enough (in particular $C$ depends on the quantity $Ce^{C V(t)}(\|q_1\|_{L^2(B^{\NN+1}_{p,1})}+\|q_2\|_{L^2(B^{\NN+1}_{p,1})})(1+C(\|q_2\|_{L^\infty})\|q_2\|_{\widetilde{L}^\infty(B^\NN_{p,1})})$) that: 
\begin{equation}
\begin{aligned}
&\|\delta q\|_{\widetilde{L}^\infty_t(B^{\NN-1}_{p,1})}+ \|\delta q\|_{\widetilde{L}^1_t(B^{\NN+1}_{p,1})}\leq C e^{C V(t)}\|\delta q\|_{\widetilde{L}^\infty_t(B^{\NN-1}_{p,1})} \int^t_0\big( \|q_2\|_{B^{\NN+2}_{p,1}}+\|{\rm div}v_2\|_{B^{\NN}_{p,1}}\\
&+\|q_1\|^2_{B^{\NN+1}_{p,1}} +( \|q_1\|^2_{B^{\NN+1}_{p,1}}+\|q_2\|^2_{B^{\NN+1}_{p,1}})(1+C(\|q_2\|_{L^\infty_t(L^\infty)})\|q_2\|_{L^\infty_t(B^\NN_{p,1})})d\tau\\
&+Ce^{C V(t)}\int^t_0 (\|{\rm div}\delta v\|_{B^{\NN-1}_{p,1}}+ \|{\rm curl}\delta v\|_{B^{\NN-1}_{p,1}})    (1+\|\bar{q}_1\|_{B^{\NN}_{p,1}}+\|\bar{q}_2\|_{B^{\NN}_{p,1}}+2\|q_L\|_{B^{\NN}_{p,1}}  )  d\tau\\
\end{aligned}
\label{transportuni2}
\end{equation}
By choosing $t$ such that:
\begin{equation}
\begin{aligned}
&C e^{C V(t)} \int^t_0\big( \|q_2\|_{B^{\NN+2}_{p,1}}+\|{\rm div}v_2\|_{B^{\NN}_{p,1}}+\|q_1\|^2_{B^{\NN+1}_{p,1}} \\
&+ (\|q_1\|^2_{B^{\NN+1}_{p,1}}+\|q_2\|^2_{B^{\NN+1}_{p,1}})(1+C(\|q_2\|_{L^\infty_t(L^\infty)})\|q_2\|_{L^\infty_t(B^\NN_{p,1})})d\tau\leq\frac{1}{2},
\end{aligned}
\label{unicondi1}
\end{equation}
we have for $C$ large enough:
\begin{equation}
\begin{aligned}
&\|\delta q\|_{\widetilde{L}^\infty_t(B^{\NN-1}_{p,1})}+ \|\delta q\|_{\widetilde{L}^1_t(B^{\NN+1}_{p,1})}\leq\\
&Ce^{C V(t)}\int^t_0 (\|{\rm div}\delta v\|_{B^{\NN-1}_{p,1}}+ \|{\rm curl}\delta v\|_{B^{\NN-1}_{p,1}})    (1+\|\bar{q}_1\|_{B^{\NN}_{p,1}}+\|\bar{q}_2\|_{B^{\NN}_{p,1}}+2\|q_L\|_{B^{\NN}_{p,1}}  )  d\tau\\
\end{aligned}
\label{transportuni3}
\end{equation}
We proceed similarly with the inequality (\ref{unicurl1}) which gives  for $C>0$ large enough: 
$$
\begin{aligned}
&\|{\rm curl}\delta v\|_{\widetilde{L}^\infty_t(B^{\NN-1}_{p,1})}+\|{\rm curl}\delta v\|_{\widetilde{L}^1_t(B^{\NN+1}_{p,1})}\leq Ce^{C V_2(t)}\int^t_0 \|\delta u\|_{B^{\NN}_{p,1}}\|{\rm curl}\bar{v}_2\|_{B^{\NN}_{p,1}}d\tau\\
&+Ce^{C V_2(t)} \|{\rm curl}\delta v\|_{\widetilde{L}_t^{\infty}(B^{\NN-1}_{p,1})}(1+\|q_1\|_{\widetilde{L}^\infty_t(B^{\NN}_{p,1})}C[\|q_1\|_{L_t^\infty(L^\infty)}])\|q_1\| _{\widetilde{L}^1_t(B^{\NN+2}_{p,1})}\\
&+Ce^{C V_2(t)}  \|{\rm curl}\delta v\|_{\widetilde{L}_t^2(B^{\NN}_{p,1})}  \|q_1\|_{L^2_t(B^{\NN+1}_{p,1})}C[\|q_1\|_{L_t^\infty(L^\infty)}]\\
&+Ce^{C V_2(t)} \|{\rm curl}\delta v\|_{\widetilde{L}^1_t(B^{\NN+1}_{p,1})} C[\|q_1\|_{L^\infty_t(L^\infty)}]   \|( Id-\dot{S}_m)q_1\|_{\widetilde{L}^\infty_t(B^\NN_{p,1})}\\
\end{aligned}
$$
\begin{equation}
\begin{aligned}
&+Ce^{C V_2(t)}\|\delta q\|_{\widetilde{L}^1_t(B^{\NN+1}_{p,1})} (1+\|q_1\|_{\widetilde{L}^\infty_t (B^{\NN}_{p,1})}C[\|q_1\|_{L^\infty}])\|{\rm curl}v_2\|_{\widetilde{L}^\infty_t(B^{\NN}_{p,1})}C[\|q_1\|_{L^\infty_t(L^\infty)}, 
\|q_2\|_{L_t^\infty(L^\infty)}]\\
&+Ce^{C V_2(t)}  \|\delta q\|_{\widetilde{L}^\infty_t(B^{\NN-1}_{p,1})} \big( \|{\rm curl}\bar{v}_2\|_{\widetilde{L}_t^1(B^{\NN+2}_{p,1})}+C[\|q_1\|_{L_t^\infty(L^\infty)}, \|q_2\|_{L^\infty_t(L^\infty)}]  \|{\rm curl}\bar{v}_2\|_{\widetilde{L}_t^1(B^{\NN+2}_{p,1})}\\
&+C[\|q_1\|_{L_t^\infty(L^\infty)}, \|q_2\|_{L_t^\infty(L^\infty)}]  \|{\rm curl}v_L\|_{\widetilde{L}_t^1(B^{\NN+2}_{p,1})}+\|{\rm curl}v_2\|_{\widetilde{L}_t^2(B^{\NN+1}_{p,1})}\|q_2\|_{\widetilde{L}_t^2(B^{\NN+1}_{p,1})}\\
&+\|q_2\|_{\widetilde{L}_t^1(B^{\NN+2}_{p,1})} C[\|q_2\|_{L_t^\infty(L^\infty)}]\|{\rm curl}v_2\|_{\widetilde{L}^\infty_t(B^{\NN}_{p,1})}\big)\\
&+Ce^{C V_2(t)} \|\delta q\|_{\widetilde{L}^2_t(B^{\NN}_{p,1})}\big( C[\|q_1\|_{L_t^\infty(L^\infty)},\|q_2\|_{L_t^\infty(L^\infty)}] \| {\rm curl}v_2\|_{\widetilde{L}^2_t(B^{\NN+1}_{p,1})}\\
&+(1+\|q_1\|_{\widetilde{L}^\infty_t(B^\NN_{p,1})}C[\|q_1\|_{L_t^\infty(L^\infty)}])\|{\rm curl}v_2\| _{\widetilde{L}^2_t(B^{\NN+1}_{p,1})}\big)\\
&+Ce^{C V_2(t)}\int^t_0 (\|\delta v\|_{B^{\NN}_{p,1}}\|u_1\|_{B^{\NN+1}_{p,1}}+ 
\|\delta u\|_{B^{\NN}_{p,1}}\|u_2\|_{B^{\NN+1}_{p,1}})d\tau\\
\end{aligned}
\label{unicurl3}
\end{equation}
Choosing $t$ small enough and $m$ large enough such that:
\begin{equation}
\begin{aligned}
&C e^{C V_2(t)}\big[(1+\|q_1\|_{\widetilde{L}^\infty_t(B^{\NN}_{p,1})}C[\|q_1\|_{L_t^\infty(L^\infty)}])\|q_1\| _{\widetilde{L}^1_t(B^{\NN+2}_{p,1})}+\|q_1\|_{L^2_t(B^{\NN+1}_{p,1})}C[\|q_1\|_{L_t^\infty(L^\infty)}]
\\
&+C[\|q_1\|_{L^\infty_t(L^\infty)}]   \|( Id-\dot{S}_m)q_1\|_{\widetilde{L}^\infty_t(B^\NN_{p,1})}]\leq\frac{1}{16},
\end{aligned}
\label{unicondi2}
\end{equation}
then by Young inequality, interpolation and $C$ large enough we have:
\begin{equation}
\begin{aligned}
&\|{\rm curl}\delta v\|_{\widetilde{L}^\infty_t(B^{\NN-1}_{p,1})}+\|{\rm curl}\delta v\|_{\widetilde{L}^1_t(B^{\NN+1}_{p,1})}\leq \\
&C\biggl( e^{C V_2(t)}\|\delta q\|_{\widetilde{L}^1_t(B^{\NN+1}_{p,1})} (1+\|q_1\|_{\widetilde{L}^\infty_t (B^{\NN}_{p,1})}C[\|q_1\|_{L^\infty}])\|{\rm curl}v_2\|_{\widetilde{L}^\infty_t(B^{\NN}_{p,1})}C[\|q_1\|_{L^\infty_t(L^\infty)}, 
\|q_2\|_{L_t^\infty(L^\infty)}]\\
&+e^{C V_2(t)}  \|\delta q\|_{\widetilde{L}^\infty_t(B^{\NN-1}_{p,1})} \big( \|{\rm curl}\bar{v}_2\|_{\widetilde{L}_t^1(B^{\NN+2}_{p,1})}+C[\|q_1\|_{L_t^\infty(L^\infty)}, \|q_2\|_{L^\infty_t(L^\infty)}]  \|{\rm curl}\bar{v}_2\|_{\widetilde{L}_t^1(B^{\NN+2}_{p,1})}\\
&+C[\|q_1\|_{L_t^\infty(L^\infty)}, \|q_2\|_{L_t^\infty(L^\infty)}]  \|{\rm curl}v_L\|_{\widetilde{L}_t^1(B^{\NN+2}_{p,1})}+\|{\rm curl}v_2\|_{\widetilde{L}_t^2(B^{\NN+1}_{p,1})}\|q_2\|_{\widetilde{L}_t^2(B^{\NN+1}_{p,1})}\\
&+\|q_2\|_{\widetilde{L}_t^1(B^{\NN+2}_{p,1})} C[\|q_2\|_{L_t^\infty(L^\infty)}]\|{\rm curl}v_2\|_{\widetilde{L}^\infty_t(B^{\NN}_{p,1})}\big)\\
&+e^{C V_2(t)} \|\delta q\|_{\widetilde{L}^2_t(B^{\NN}_{p,1})}\big( C[\|q_1\|_{L_t^\infty(L^\infty)},\|q_2\|_{L_t^\infty(L^\infty)}] \| {\rm curl}v_2\|_{\widetilde{L}^2_t(B^{\NN+1}_{p,1})}\\
&+(1+\|q_1\|_{\widetilde{L}^\infty_t(B^\NN_{p,1})}C[\|q_1\|_{L_t^\infty(L^\infty)}])\|{\rm curl}v_2\| _{\widetilde{L}^2_t(B^{\NN+1}_{p,1})})\biggl)\\
&+C e^{C V_2(t)}\int^t_0 (\|\delta v\|_{B^{\NN}_{p,1}}\|u_1\|_{B^{\NN+1}_{p,1}}+ 
\|\delta u\|_{B^{\NN}_{p,1}}\|u_2\|_{B^{\NN+1}_{p,1}})d\tau\big).
\end{aligned}
\label{unicurl4}
\end{equation}
We recall now by definition of $v$ (since $v=u+\n\va(\rho)$) and proposition \ref{estimcompo} that for $C$ large enough:
\begin{equation}
\begin{aligned}
&\|\delta v\|_{B^{\NN}_{p,1}}\leq C(\|{\rm div}\delta v\|_{B^{\NN-1}_{p,1}}+\|{\rm curl}\delta v\|_{B^{\NN-1}_{p,1}}),\\
&\|\delta u\|_{B^{\NN}_{p,1}}\leq C(\|{\rm div}\delta v\|_{B^{\NN-1}_{p,1}}+\|{\rm curl}\delta v\|_{B^{\NN-1}_{p,1}}+C[\|q_1\|_{L^\infty}, \|q_2\|_{L^\infty}]\|\delta q\|_{B^{\NN+1}_{p,1}}),
\end{aligned}
\label{diffuv}
\end{equation}
putting this inequality in (\ref{unicurl4}) we have:
\begin{equation}
\begin{aligned}
&\|{\rm curl}\delta v\|_{\widetilde{L}^\infty_t(B^{\NN-1}_{p,1})}+\|{\rm curl}\delta v\|_{\widetilde{L}^1_t(B^{\NN+1}_{p,1})}\leq \\
&C\biggl( e^{C V_2(t)}\|\delta q\|_{\widetilde{L}^1_t(B^{\NN+1}_{p,1})} (1+\|q_1\|_{\widetilde{L}^\infty_t (B^{\NN}_{p,1})}C[\|q_1\|_{L^\infty}])\|{\rm curl}v_2\|_{\widetilde{L}^\infty_t(B^{\NN}_{p,1})}C[\|q_1\|_{L^\infty_t(L^\infty)}, 
\|q_2\|_{L_t^\infty(L^\infty)}]\\
&+e^{C V_2(t)}  \|\delta q\|_{\widetilde{L}^\infty_t(B^{\NN-1}_{p,1})} \big( \|{\rm curl}\bar{v}_2\|_{\widetilde{L}_t^1(B^{\NN+2}_{p,1})}+C[\|q_1\|_{L_t^\infty(L^\infty)}, \|q_2\|_{L^\infty_t(L^\infty)}]  \|{\rm curl}\bar{v}_2\|_{\widetilde{L}_t^1(B^{\NN+2}_{p,1})}\\
&+C[\|q_1\|_{L_t^\infty(L^\infty)}, \|q_2\|_{L_t^\infty(L^\infty)}]  \|{\rm curl}v_L\|_{\widetilde{L}_t^1(B^{\NN+2}_{p,1})}+\|{\rm curl}v_2\|_{\widetilde{L}_t^2(B^{\NN+1}_{p,1})}\|q_2\|_{\widetilde{L}_t^2(B^{\NN+1}_{p,1})}\\
&+\|q_2\|_{\widetilde{L}_t^1(B^{\NN+2}_{p,1})} C[\|q_2\|_{L_t^\infty(L^\infty)}]\|{\rm curl}v_2\|_{\widetilde{L}^\infty_t(B^{\NN}_{p,1})}\big)\\
&+e^{C V_2(t)} \|\delta q\|_{\widetilde{L}^2_t(B^{\NN}_{p,1})}\big( C[\|q_1\|_{L_t^\infty(L^\infty)},\|q_2\|_{L_t^\infty(L^\infty)}] \| {\rm curl}v_2\|_{\widetilde{L}^2_t(B^{\NN+1}_{p,1})}\\
&+(1+\|q_1\|_{\widetilde{L}^\infty_t(B^\NN_{p,1})}C[\|q_1\|_{L_t^\infty(L^\infty)}])\|{\rm curl}v_2\| _{\widetilde{L}^2_t(B^{\NN+1}_{p,1})})\biggl)\\
&+C e^{C V_2(t)}\int^t_0 \big((\|{\rm div}\delta v\|_{B^{\NN-1}_{p,1}}+\|{\rm curl}\delta v\|_{B^{\NN-1}_{p,1}})
(\|u_1\|_{B^{\NN+1}_{p,1}}+ \|u_2\|_{B^{\NN+1}_{p,1}})\\
&\hspace{8cm}+\|\delta q\|_{B^{\NN+1}_{p,1}} \|u_2\|_{B^{\NN+1}_{p,1}}     \big)d\tau\big).
\end{aligned}
\label{unicurl5}
\end{equation}
Assuming $t$ small enough such that:
\begin{equation}
C e^{C V_2(t)}(\|u_1\|_{\widetilde{L}^1_t(B^{\NN+1}_{p,1})}+ \|u_2\|_{\widetilde{L}^1_t(B^{\NN+1}_{p,1})})\leq\frac{1}{16},
\label{unicondi6}
\end{equation}
we have for $C$ large enough:
\begin{equation}
\begin{aligned}
&\|{\rm curl}\delta v\|_{\widetilde{L}^\infty_t(B^{\NN-1}_{p,1})}+\|{\rm curl}\delta v\|_{\widetilde{L}^1_t(B^{\NN+1}_{p,1})}\leq \\
&C\big( e^{C V_2(t)}\|\delta q\|_{\widetilde{L}^1_t(B^{\NN+1}_{p,1})} (1+\|q_1\|_{\widetilde{L}^\infty_t (B^{\NN}_{p,1})}C[\|q_1\|_{L^\infty}])\|{\rm curl}v_2\|_{\widetilde{L}^\infty_t(B^{\NN}_{p,1})}C[\|q_1\|_{L^\infty_t(L^\infty)}, 
\|q_2\|_{L_t^\infty(L^\infty)}]\\
&+e^{C V_2(t)}  \|\delta q\|_{\widetilde{L}^\infty_t(B^{\NN-1}_{p,1})} \big( \|{\rm curl}\bar{v}_2\|_{\widetilde{L}_t^1(B^{\NN+2}_{p,1})}+C[\|q_1\|_{L_t^\infty(L^\infty)}, \|q_2\|_{L^\infty_t(L^\infty)}]  \|{\rm curl}\bar{v}_2\|_{\widetilde{L}_t^1(B^{\NN+2}_{p,1})}\\
&+C[\|q_1\|_{L_t^\infty(L^\infty)}, \|q_2\|_{L_t^\infty(L^\infty)}]  \|{\rm curl}v_L\|_{\widetilde{L}_t^1(B^{\NN+2}_{p,1})}+\|{\rm curl}v_2\|_{\widetilde{L}_t^2(B^{\NN+1}_{p,1})}\|q_2\|_{\widetilde{L}_t^2(B^{\NN+1}_{p,1})}\\
&+\|q_2\|_{\widetilde{L}_t^1(B^{\NN+2}_{p,1})} C[\|q_2\|_{L_t^\infty(L^\infty)}]\|{\rm curl}v_2\|_{\widetilde{L}^\infty_t(B^{\NN}_{p,1})}\big)\\
&+e^{C V_2(t)} \|\delta q\|_{\widetilde{L}^2_t(B^{\NN}_{p,1})}\big( C[\|q_1\|_{L_t^\infty(L^\infty)},\|q_2\|_{L_t^\infty(L^\infty)}] \| {\rm curl}v_2\|_{\widetilde{L}^2_t(B^{\NN+1}_{p,1})}\\
&+(1+\|q_1\|_{\widetilde{L}^\infty_t(B^\NN_{p,1})}C[\|q_1\|_{L_t^\infty(L^\infty)}])\|{\rm curl}v_2\| _{\widetilde{L}^2_t(B^{\NN+1}_{p,1})}\big)\\
&+e^{C V_2(t)}\int^t_0 \big(\|{\rm div}\delta v\|_{B^{\NN-1}_{p,1}}
(\|u_1\|_{B^{\NN+1}_{p,1}}+ \|u_2\|_{B^{\NN+1}_{p,1}})+\|\delta q\|_{B^{\NN+1}_{p,1}} \|u_2\|_{B^{\NN+1}_{p,1}}     \big)d\tau\big).
\end{aligned}
\label{unicurl6}
\end{equation}
Using (\ref{transportuni3}) and plugging in (\ref{unicurl6}) we have for $C>0$ large enough:
\begin{equation}
\begin{aligned}
&\|{\rm curl}\delta v\|_{\widetilde{L}^\infty_t(B^{\NN-1}_{p,1})}+\|{\rm curl}\delta v\|_{\widetilde{L}^1_t(B^{\NN+1}_{p,1})}\leq \\
&C e^{C V_2(t)}e^{C V(t)}\int^t_0 (\|{\rm div}\delta v\|_{B^{\NN-1}_{p,1}}+ \|{\rm curl}\delta v\|_{B^{\NN-1}_{p,1}})    (1+\|\bar{q}_1\|_{B^{\NN}_{p,1}}+\|\bar{q}_2\|_{B^{\NN}_{p,1}}+2\|q_L\|_{B^{\NN}_{p,1}}  )  d\tau\\
&\times\biggl[ (1+\|q_1\|_{\widetilde{L}^\infty_t (B^{\NN}_{p,1})}C[\|q_1\|_{L^\infty}])\|{\rm curl}v_2\|_{\widetilde{L}^\infty_t(B^{\NN}_{p,1})}C[\|q_1\|_{L^\infty_t(L^\infty)}, 
\|q_2\|_{L_t^\infty(L^\infty)}]\\
&+ \big( \|{\rm curl}\bar{v}_2\|_{\widetilde{L}_t^1(B^{\NN+2}_{p,1})}+C[\|q_1\|_{L_t^\infty(L^\infty)}, \|q_2\|_{L^\infty_t(L^\infty)}]  \|{\rm curl}\bar{v}_2\|_{\widetilde{L}_t^1(B^{\NN+2}_{p,1})}\\
&+C[\|q_1\|_{L_t^\infty(L^\infty)}, \|q_2\|_{L_t^\infty(L^\infty)}]  \|{\rm curl}v_L\|_{\widetilde{L}_t^1(B^{\NN+2}_{p,1})}+\|{\rm curl}v_2\|_{\widetilde{L}_t^2(B^{\NN+1}_{p,1})}\|q_2\|_{\widetilde{L}_t^2(B^{\NN+1}_{p,1})}\\
&+\|q_2\|_{\widetilde{L}_t^1(B^{\NN+2}_{p,1})} C[\|q_2\|_{L_t^\infty(L^\infty)}]\|{\rm curl}v_2\|_{\widetilde{L}^\infty_t(B^{\NN}_{p,1})}\big)\\
&+\big( C[\|q_1\|_{L_t^\infty(L^\infty)},\|q_2\|_{L_t^\infty(L^\infty)}] \| {\rm curl}v_2\|_{\widetilde{L}^2_t(B^{\NN+1}_{p,1})}\\
&+(1+\|q_1\|_{\widetilde{L}^\infty_t(B^\NN_{p,1})}C[\|q_1\|_{L_t^\infty(L^\infty)}])\|{\rm curl}v_2\| _{\widetilde{L}^2_t(B^{\NN+1}_{p,1})}\big)+  \|u_2\|_{\widetilde{L}^\infty_t(B^{\NN+1}_{p,1})}\biggl]\\
&+e^{C V_2(t)}\int^t_0 \big(\|{\rm div}\delta v\|_{B^{\NN-1}_{p,1}}
(\|u_1\|_{B^{\NN+1}_{p,1}}+ \|u_2\|_{B^{\NN+1}_{p,1}})    \big)d\tau\big).
\end{aligned}
\label{unicurl7}
\end{equation}
In the same way than previously, choosing $t$ small enough such that:
\begin{equation}
\begin{aligned}
&C e^{C V_2(t)}e^{C V(t)}\int^t_0 (1+\|\bar{q}_1\|_{B^{\NN}_{p,1}}+\|\bar{q}_2\|_{B^{\NN}_{p,1}}+2\|q_L\|_{B^{\NN}_{p,1}}  )  d\tau\\
&\times\biggl[ (1+\|q_1\|_{\widetilde{L}^\infty_t (B^{\NN}_{p,1})}C[\|q_1\|_{L^\infty}])\|{\rm curl}v_2\|_{\widetilde{L}^\infty_t(B^{\NN}_{p,1})}C[\|q_1\|_{L^\infty_t(L^\infty)}, 
\|q_2\|_{L_t^\infty(L^\infty)}]\\
&+ \big( \|{\rm curl}\bar{v}_2\|_{\widetilde{L}_t^1(B^{\NN+2}_{p,1})}+C[\|q_1\|_{L_t^\infty(L^\infty)}, \|q_2\|_{L^\infty_t(L^\infty)}]  \|{\rm curl}\bar{v}_2\|_{\widetilde{L}_t^1(B^{\NN+2}_{p,1})}\\
&+C[\|q_1\|_{L_t^\infty(L^\infty)}, \|q_2\|_{L_t^\infty(L^\infty)}]  \|{\rm curl}v_L\|_{\widetilde{L}_t^1(B^{\NN+2}_{p,1})}+\|{\rm curl}v_2\|_{\widetilde{L}_t^2(B^{\NN+1}_{p,1})}\|q_2\|_{\widetilde{L}_t^2(B^{\NN+1}_{p,1})}\\
&+\|q_2\|_{\widetilde{L}_t^1(B^{\NN+2}_{p,1})} C[\|q_2\|_{L_t^\infty(L^\infty)}]\|{\rm curl}v_2\|_{\widetilde{L}^\infty_t(B^{\NN}_{p,1})}\big)\\
&+\big( C[\|q_1\|_{L_t^\infty(L^\infty)},\|q_2\|_{L_t^\infty(L^\infty)}] \| {\rm curl}v_2\|_{\widetilde{L}^2_t(B^{\NN+1}_{p,1})}\\
&+(1+\|q_1\|_{\widetilde{L}^\infty_t(B^\NN_{p,1})}C[\|q_1\|_{L_t^\infty(L^\infty)}])\|{\rm curl}v_2\| _{\widetilde{L}^2_t(B^{\NN+1}_{p,1})}\big)+  \|u_2\|_{\widetilde{L}^\infty_t(B^{\NN+1}_{p,1})}\biggl]\leq\frac{1}{16},
\end{aligned}
\label{unicondi3}
\end{equation}
we have for $C>0$ large enough:
\begin{equation}
\begin{aligned}
&\|{\rm curl}\delta v\|_{\widetilde{L}^\infty_t(B^{\NN-1}_{p,1})}+\|{\rm curl}\delta v\|_{\widetilde{L}^1_t(B^{\NN+1}_{p,1})}\leq \\
&C e^{C V_2(t)}e^{C V(t)}\int^t_0 \|{\rm div}\delta v\|_{B^{\NN-1}_{p,1}}   \big(1+\|\bar{q}_1\|_{B^{\NN}_{p,1}}+\|\bar{q}_2\|_{B^{\NN}_{p,1}}+2\|q_L\|_{B^{\NN}_{p,1}}\\
&\hspace{7cm}+\|u_1\|_{B^{\NN+1}_{p,1}}+ \|u_2\|_{B^{\NN+1}_{p,1}}  \big)  d\tau\\
&\times\biggl[ (1+\|q_1\|_{\widetilde{L}^\infty_t (B^{\NN}_{p,1})}C[\|q_1\|_{L^\infty}])\|{\rm curl}v_2\|_{\widetilde{L}^\infty_t(B^{\NN}_{p,1})}C[\|q_1\|_{L^\infty_t(L^\infty)}, 
\|q_2\|_{L_t^\infty(L^\infty)}]\\
&+ \big( \|{\rm curl}\bar{v}_2\|_{\widetilde{L}_t^1(B^{\NN+2}_{p,1})}+C[\|q_1\|_{L_t^\infty(L^\infty)}, \|q_2\|_{L^\infty_t(L^\infty)}]  \|{\rm curl}\bar{v}_2\|_{\widetilde{L}_t^1(B^{\NN+2}_{p,1})}\\
&+C[\|q_1\|_{L_t^\infty(L^\infty)}, \|q_2\|_{L_t^\infty(L^\infty)}]  \|{\rm curl}v_L\|_{\widetilde{L}_t^1(B^{\NN+2}_{p,1})}+\|{\rm curl}v_2\|_{\widetilde{L}_t^2(B^{\NN+1}_{p,1})}\|q_2\|_{\widetilde{L}_t^2(B^{\NN+1}_{p,1})}\\
&+\|q_2\|_{\widetilde{L}_t^1(B^{\NN+2}_{p,1})} C[\|q_2\|_{L_t^\infty(L^\infty)}]\|{\rm curl}v_2\|_{\widetilde{L}^\infty_t(B^{\NN}_{p,1})}\big)\\
&+\big( C[\|q_1\|_{L_t^\infty(L^\infty)},\|q_2\|_{L_t^\infty(L^\infty)}] \| {\rm curl}v_2\|_{\widetilde{L}^2_t(B^{\NN+1}_{p,1})}\\
&+(1+\|q_1\|_{\widetilde{L}^\infty_t(B^\NN_{p,1})}C[\|q_1\|_{L_t^\infty(L^\infty)}])\|{\rm curl}v_2\| _{\widetilde{L}^2_t(B^{\NN+1}_{p,1})}\big)+  \|u_2\|_{\widetilde{L}^\infty_t(B^{\NN+1}_{p,1})}\biggl].
\end{aligned}
\label{unicurl8}
\end{equation}
We plug now in (\ref{unidiv}) (\ref{unicurl8}) and (\ref{transportuni1}), it yields:
$$
\begin{aligned}
&\|{\rm div}\delta v\|_{\widetilde{L}^\infty_t(B^{\NN-1}_{p,1})}\leq C e^{C V_1(t)}\int^t_0\big(\|\delta v\|_{B^{\NN}_{p,1}} \|u_1\|_{B^{\NN+1}_{p,1}} +  \|\delta u\|_{B^\NN_{p,1}}(\|{\rm div}\bar{v}_2\|_{B^{\NN}_{p,1}}+ \|v_2\|_{B^{\NN+1}_{p,1}})\big) d\tau\\
&e^{C(V(t)+V_1(t))  }   \int^t_0(\|{\rm div}\delta v\|_{B^{\NN-1}_{p,1}}+ \|{\rm curl}\delta v\|_{B^{\NN-1}_{p,1}})(1+\|\bar{q}_1\|_{B^{\NN}_{p,1}}+\|\bar{q}_2\|_{B^{\NN}_{p,1}}+2\|q_L\|_{B^{\NN}_{p,1}}  )  d\tau
\\
&\times\biggl[\|q_2\|_{\widetilde{L}^2_t(B^{\NN+1}_{p,1})}   \|{\rm curl}v_2 \|_{\widetilde{L}^2_t(B^{\NN+1}_{p,1})} \\
&+( \|{\rm div}{\rm curl}v_1\|_{\widetilde{L}^2_t(B^{\NN}_{p,1})}\|q_1\|_{\widetilde{L}^\infty_t(B^\NN_{p,1})}C[\|q\|_{L_t^\infty(L^\infty)}] + \| {\rm curl}v_2\|_{\widetilde{L}^2_t(B^{\NN}_{p,1})})\\
&+\big( C[\|q_1\|_{L_t^\infty(L^\infty)}, \|q_2\|_{L_t^\infty(L^\infty)}] +C[\|q_1\|_{L_t^\infty(L^\infty)}, \|q_2\|_{L_t^\infty(L^\infty)}] \|{\rm curl}v_1\|_{\widetilde{L}^\infty_t(B^{\NN-1}_{p,1})})\biggl]\\[4mm]
&+Ce^{C V_1(t)} e^{C V_2(t)}e^{C V(t)}\int^t_0 \|{\rm div}\delta v\|_{B^{\NN-1}_{p,1}}   \big(1+\|\bar{q}_1\|_{B^{\NN}_{p,1}}+\|\bar{q}_2\|_{B^{\NN}_{p,1}}+2\|q_L\|_{B^{\NN}_{p,1}}\\
&\hspace{9cm}+\|u_1\|_{B^{\NN+1}_{p,1}}+ \|u_2\|_{B^{\NN+1}_{p,1}}  \big)  d\tau\\
&\times\biggl[ (1+\|q_1\|_{\widetilde{L}^\infty_t (B^{\NN}_{p,1})}C[\|q_1\|_{L^\infty}])\|{\rm curl}v_2\|_{\widetilde{L}^\infty_t(B^{\NN}_{p,1})}C[\|q_1\|_{L^\infty_t(L^\infty)}, 
\|q_2\|_{L_t^\infty(L^\infty)}]\\
&+ \big( \|{\rm curl}\bar{v}_2\|_{\widetilde{L}_t^1(B^{\NN+2}_{p,1})}+C[\|q_1\|_{L_t^\infty(L^\infty)}, \|q_2\|_{L^\infty_t(L^\infty)}]  \|{\rm curl}\bar{v}_2\|_{\widetilde{L}_t^1(B^{\NN+2}_{p,1})}
\end{aligned}
$$
\begin{equation}
\begin{aligned}
&+C[\|q_1\|_{L_t^\infty(L^\infty)}, \|q_2\|_{L_t^\infty(L^\infty)}]  \|{\rm curl}v_L\|_{\widetilde{L}_t^1(B^{\NN+2}_{p,1})}+\|{\rm curl}v_2\|_{\widetilde{L}_t^2(B^{\NN+1}_{p,1})}\|q_2\|_{\widetilde{L}_t^2(B^{\NN+1}_{p,1})}\\
&+\|q_2\|_{\widetilde{L}_t^1(B^{\NN+2}_{p,1})} C[\|q_2\|_{L_t^\infty(L^\infty)}]\|{\rm curl}v_2\|_{\widetilde{L}^\infty_t(B^{\NN}_{p,1})}\big)\\
&+\big( C[\|q_1\|_{L_t^\infty(L^\infty)},\|q_2\|_{L_t^\infty(L^\infty)}] \| {\rm curl}v_2\|_{\widetilde{L}^2_t(B^{\NN+1}_{p,1})}\\
&+(1+\|q_1\|_{\widetilde{L}^\infty_t(B^\NN_{p,1})}C[\|q_1\|_{L_t^\infty(L^\infty)}])\|{\rm curl}v_2\| _{\widetilde{L}^2_t(B^{\NN+1}_{p,1})}\big)+  \|u_2\|_{\widetilde{L}^\infty_t(B^{\NN+1}_{p,1})}\biggl]\\
&\times\biggl[\|q_2\|_{\widetilde{L}^2_t(B^{\NN+1}_{p,1})}+\|q_1\|_{\widetilde{L}^2_t(B^{\NN}_{p,1})} C[\|q\|_{L_t^\infty(L^\infty)}])+ C[\|q_2\|_{L_t^\infty(L^\infty)}] \|q_2\|_{\widetilde{L}^1_t(B^{\NN+2}_{p,1})}
\biggl]\\
&
\end{aligned}
\label{unidiv2}
\end{equation}
Plugging once more the inequalities  (\ref{diffuv}), (\ref{unicurl8}), (\ref{transportuni3}) in (\ref{unidiv2}), and using the Gronwall inequality we conclude that:
$$\|{\rm div}\delta v\|_{\widetilde{L}^\infty_t(B^{\NN-1}_{p,1})}=0\;\;\mbox{for any}\;t\in[0,T],  $$
with $T$ verifying (\ref{unicondi1}), (\ref{unicondi2}), (\ref{unicondi6}) and (\ref{unicondi3}). We deduce using (\ref{transportuni3}) (\ref{unicurl8})  that:
$$\|{\rm curl}\delta v\|_{\widetilde{L}^\infty_t(B^{\NN-1}_{p,1})}=\|\delta q\|_{\widetilde{L}^\infty_t(B^{\NN}_{p,1})} =0\;\;\mbox{for any}\;t\in[0,T].$$
It implies that $q_1=q_2$ and $v_1=v_2$ on $[0,T]$.
\\
To end the proof when $T$ is not small, let us introduce (as in \cite{DW}, section $10.2.4$) the set:
$$
I\overset{\mbox{def}}{=}\{t\in[0,T]/(q_1(t'), \ub_1(t'))=(q_2(t'), \ub_2(t')),\;\forall t'\in[0,t]\}.
$$
This is a nonempty closed subset of $[0,T]$. Using the same method as above allows to prove it is also open and then $I=[0,T]$.

\section{Global well-posedness, Proof of the theorem \ref{theo2}}
\label{section5}
In this section we are interested in proving the global well-posedness of (\ref{0.1}) when we assume smallness on the initial data. The proof follows the same lines than in the sections \ref{section3} and \ref{section4}. The main difficulty consists in getting $v$ and $u$ in $\widetilde{L}^1(\R^{+},B^{\NN+1}_{p,1})$ or at least to prove that $\n u$ and $\n v$ belong in $L^1(\R^+, L^\infty)$ (indeed we need to take into account the behavior in low frequencies). This is necessary since ${\rm div}v$ verifies a transport equation, and if we want to propagate the regularity on ${\rm div}v$, it implies a control on the Lipschitz norm of $u$. In the previous section we have only an estimate of ${\rm div}v$ in $\widetilde{L}_T^\infty(B^{\NN}_{p,1})$ which is not sufficient for large time, indeed we are interested in getting an $L^1$ estimate in time. It is necessary to exhibit a damped effect on ${\rm div}v$. In addition in the previous section ${\rm curl}v$ was only in $\widetilde{L}_T^1(B^{\NN+2}_{p,1})\cap \widetilde{L}^\infty(B^{\NN}_{p,1})$, for the same reason we need to control  ${\rm curl}v$ in $\widetilde{L}^1(\R^+,B^{\NN}_{p,1})$. It requires then additional assumption in low frequencies on ${\rm curl}v_0$. More precisely we must assume that ${\rm curl}v_0$ is in $B^{\N-2}_{2,1}$ in low frequencies (indeed we have to work with $p=2$ since in low frequencies the system has a hyperbolic behavior) in order to ensure a control of  ${\rm curl}v$ in $\widetilde{L}^1(B^{\N}_{2,1})$. All these considerations explain the choice on the initial data.\\
In particular we observe that the conditions on the initial data of theorem \ref{theo2} are the same than in \cite{arma}, it implies that we have the existence of global strong solution and in addition we have the regularity that we wish on $(q,{\rm div}v,{\rm curl}v)$ in low frequencies. The only thing which remains to do is to prove the regularity assumptions on $q$, ${\rm div}v$ and ${\rm curl}v$ in high frequencies. To do this we have just to  estimates in Besov spaces the following linear system associated to (\ref{0.6}) in high frequencies:
\begin{equation}
\begin{cases}
\begin{aligned}
&\p_{t}q-2\mu\D q+{\rm div}v=F,\\
&\p_{t}{\rm div}v+u\cdot\n {\rm div}v+\D q=G,\\
&\p_t{\rm curl}v-\mu\D{\rm curl}v=H,\\
\end{aligned}
\end{cases}
\label{lineaire}
\end{equation}
Let us study this system and let us give the following proposition which is inspired from 
 \cite{CD,arma} .
\begin{proposition}
\label{5linear1} Let $(q,v)$ a solution of the system (\ref{lineaire}) on
$[0,T[$  and $1\leq p<\max(4,N)$,
$V(t)=\int^{t}_{0}\|\nabla u(\tau)\|_{\widetilde{B}^{\N+1,\NN+1}_{2,p,1}}d\tau$. We have then the
following estimate for any $T>0$:
$$
\begin{aligned}
&\|(q,v)\|_{\widetilde{L}^{\infty}_{T}(\widetilde{B}^{\N-1,\NN}_{2,p,1})\times \widetilde{L}^{\infty}_{T}(\widetilde{B}^{\N-1,\NN+1}_{2,p,1})}+\|(q,{\rm curl}v,{\rm div}v)\|_{\widetilde{L}^{1}_{T}(\widetilde{B}^{\N+1,\NN+2}_{2,p,1})\times \widetilde{L}^{1}_{T}(\widetilde{B}^{\N,\NN+2}_{2,p,1}) \times \widetilde{L}^{1}_{T}(\widetilde{B}^{\N,\NN}_{2,p,1})}\\
&\hspace{1cm}\leq Ce^{CV(t)}\big(\|
(q_{0},v_{0})\|_{\widetilde{B}^{\N-1,\NN}_{2,p,1}\times \widetilde{B}^{\N-1,\NN+1}_{2,p,1} }\\
&\hspace{2cm}+\int^{T}_{0}
e^{-CV(\tau)}\|
(F,G,H)(\tau)\|_{\widetilde{B}^{\N-1,\NN}_{2,p,1}\times \widetilde{B}^{\N-2,\NN}_{2,p,1}\times \widetilde{B}^{\N-2,\NN}_{2,p,1}}d\tau\big),\\
\end{aligned}
$$
where $C$ depends only on  $N$.
\end{proposition}
{\bf Proof:} The third equation from (\ref{lineaire}) is just a heat equation, it suffices then to apply the proposition \ref{chaleur}. Let us deal with the two first equation, applying the gradient to the first equation we have:
\begin{equation}
\begin{cases}
\begin{aligned}
&\p_{t}\n q-2\mu\D \n q+\n {\rm div}v=\n F,\\
&\p_{t}{\rm div}v+u\cdot\n {\rm div}v+{\rm div}\n q=G.
\end{aligned}
\end{cases}
\label{lineaire1}
\end{equation}
Setting $q'={\rm div}v$ and $u'=\n q$, we have to study the system:
\begin{equation}
\begin{cases}
\begin{aligned}
&\p_{t}q'+u\cdot\n q'+{\rm div}u'=G,\\
&\p_{t}u'-2\mu\D u'+\n q'=\n F.
\end{aligned}
\end{cases}
\label{lineaire1}
\end{equation}
This system has been studied in \cite{CD,arma} and it is proven that for any $T>0$:
\begin{equation}
\begin{aligned}
&\|q'\|_{\widetilde{L}^\infty(\widetilde{B}^{\N-1,\NN}_{2,p,1})}+\|q'\|_{\widetilde{L}^1(\widetilde{B}^{\N+1,\NN}_{2,p,1})}+\|u'\|_{\widetilde{L}^\infty(\widetilde{B}^{\N-1,\NN-1}_{2,p,1})}+\|u'\|_{\widetilde{L}^1(\widetilde{B}^{\N+1,\NN+1}_{2,p,1})}\\
&\leq C e^{V(T)}(\|q'_0\|_{\widetilde{B}^{\N-1,\NN}_{2,p,1}}+\|u'_0\|_{\widetilde{B}^{\N-1,\NN-1}_{2,p,1}}+\|G\|_{\widetilde{L}^1(\widetilde{B}^{\N-1,\NN}_{2,p,1})}+\| \n F\| _{\widetilde{L}^1(\widetilde{B}^{\N-1,\NN-1}_{2,p,1})}   ).
\end{aligned}
\label{aestim1}
\end{equation}
In other words it gives:
\begin{equation}
\begin{aligned}
&\|{\rm div}v\|_{\widetilde{L}^\infty(\widetilde{B}^{\N-1,\NN}_{2,p,1})}+\|{\rm div}v\|_{\widetilde{L}^1(\widetilde{B}^{\N+1,\NN}_{2,p,1})}+\|q\|_{\widetilde{L}^\infty(\widetilde{B}^{\N,\NN}_{2,p,1})}+\|q\|_{\widetilde{L}^1(\widetilde{B}^{\N+2,\NN+2}_{2,p,1})}\\
&\leq C e^{V(T)}(\|{\rm div}v_0\|_{\widetilde{B}^{\N-1,\NN}_{2,p,1}}+\|q_0\|_{\widetilde{B}^{\N,\NN}_{2,p,1}}+  \|G\|_{\widetilde{L}^1(\widetilde{B}^{\N-1,\NN}_{2,p,1})}+\| \n F\| _{\widetilde{L}^1(\widetilde{B}^{\N-1,\NN-1}_{2,p,1})} ).
\end{aligned}
\label{aestim2}
\end{equation}
In a similar way we have:
\begin{equation}
\begin{aligned}
&\|{\rm div}v\|_{\widetilde{L}^\infty(\widetilde{B}^{\N-2,\NN-1}_{2,p,1})}+\|{\rm div}v\|_{\widetilde{L}^1(\widetilde{B}^{\N,\NN-1}_{2,p,1})}+\|q\|_{\widetilde{L}^\infty(\widetilde{B}^{\N-1,\NN-1}_{2,p,1})}+\|q\|_{\widetilde{L}^1(\widetilde{B}^{\N+1,\NN+1}_{2,p,1})}\\
&\leq C e^{V(T)}(\|{\rm div}v_0\|_{\widetilde{B}^{\N-2,\NN-1}_{2,p,1}}+\|q_0\|_{\widetilde{B}^{\N-1,\NN-1}_{2,p,1}}+  \|G\|_{\widetilde{L}^1(\widetilde{B}^{\N-2,\NN-1}_{2,p,1})}+\| \n F\| _{\widetilde{L}^1(\widetilde{B}^{\N-2,\NN-2}_{2,p,1})} ).
\end{aligned}
\label{aestim2}
\end{equation}Collecting (\ref{aestim1}), (\ref{aestim2}) and the estimates on ${\rm curl}v$ via the heat equation, we have the desired result.  $\blacksquare$\\
\\
The rest of the proof consists in searching a solution of the form $(q,{\rm div}v,{\rm curl}v)$ and using the proposition \ref{5linear1} on $(q,{\rm div}v,{\rm curl}v)$ and to bound the remainder term $(F,G,H)$ in $\widetilde{L}^1(\R^+,\widetilde{B}^{\N-1,\NN}_{2,p,1}\times \widetilde{B}^{\N-2,\NN}_{2,p,1}\times \widetilde{B}^{\N-2,\NN}_{2,p,1})$.
We should to treat the remainder terms as in section \ref{section3}. For the uniqueness the method follows the same approach as in section \ref{section4}.
\section{Appendix}
\label{appendix}
Here we just start by explaining how we pass from the system (\ref{0.1ab}) to the system (\ref{1.b}), we are going by giving some details on the computations.  Indeed we apply the operators ${\rm div}$ and ${\rm curl}$ to the momentum equation of (\ref{0.1a}), we have in particular:
$$
\begin{aligned}
&{\rm div}(u\cdot\n v)= u\cdot\n {\rm div}v+\n v: ^{t}\n u,\\
&{\rm curl}(u\cdot\n v)_{ij}=u\cdot\n({\rm curl}v)_{ij}+\sum_k (\p_i u_k\p_k v_j-\p_j u_k\p_k v_i),\\
&{\rm curl}(\frac{\mu(\rho)}{\rho}{\rm div}{\rm curl}v)_{ij}=\frac{\mu(\rho)}{\rho}\D {\rm curl}v_{ij}+\p_i(\frac{\mu(\rho)}{\rho})\sum_k \p_k({\rm curl}_{kj})-\p_j(\frac{\mu(\rho)}{\rho})\sum_k \p_k({\rm curl}_{ki})\\
&=\frac{\mu(\rho)}{\rho}\D {\rm curl}v_{ij}+\frac{1}{2}\frac{\lambda(\rho)}{\rho^2}\big(\sum_k (\p_i\rho\p_k({\rm curl}_{kj})-\p_j\rho\p_k({\rm curl}_{ki})\big),\\[2mm]
&({\rm curl}(\n\va(\rho)\cdot{\rm curl}v))_{kj}=\n\va(\rho)\cdot ({\rm curl}v)_{kj}+\sum_i \big(\p_{ik}\va(\rho)({\rm curl}v)_{ij}-\p_{ij}\va(\rho)({\rm curl}v)_{ik}\big).
\end{aligned}
$$
These calculus allow to obtain the system (\ref{1.b}). Let us explain now why the system (\ref{0.1}) is equivalent to the system (\ref{0.1a}) when we introduce the effective velocity $v=u+\n\va(\rho)$.
\begin{proposition}
We can formally rewrite the system (\ref{0.1}) when $\mu$ and $\lambda$ verify the condition (\ref{1.2}) as follows with $v=u+\n\va(\rho)$ ($\va'(\rho)=\frac{2\mu'(\rho)}{\rho}$):
\begin{equation}
\begin{cases}
\begin{aligned}
&\p_{t}\rho-2\D\mu(\rho)+{\rm div}(\rho v)=0,\\
&\rho \p_{t}v+\rho u\cdot\n v-{\rm div}(\mu(\rho){\rm curl}v)+\n P(\rho)=0,\\
&(\rho,u)_{/t=0}=(\rho_{0},u_{0}).
\end{aligned}
\end{cases}
\label{0.1a}
\end{equation}
\end{proposition}
{\bf Proof:}
As observed in \cite{cras}, we are interested in rewriting the system (\ref{0.1}) in terms of the following unknown $v=u+\n f(\rho)$ where $f$ is a regular function. We have then by using the mass equation:
$$
\p_t\n f(\rho)+\n\big( f'(\rho)\rho\,{\rm div}u+u\cdot\n f(\rho)\big)=0,$$
and we have:
$$\rho\p_t\n f(\rho)+\rho\n(f'(\rho)\rho{\rm div}u)+\rho u\cdot\n\n f(\rho)+\rho\n f(\rho)\cdot ^t\n u=0.
$$
By summing the previous equality with the momentum equation, we get by setting $v=u+\n f(\rho)$:
$$
\begin{aligned}
&\rho\p_t v+\rho u\cdot\n v-\mu(\rho)\D u-\mu(\rho)\n{\rm div}u-2\n\mu(\rho)\cdot Du-\lambda(\rho)\n{\rm div}u\\
&\hspace{3cm}-\n\lambda(\rho)\, {\rm div}u+\rho\n(f'(\rho)\rho{\rm div}u)+\rho\n f(\rho)\cdot ^t\n u
+\n P(\rho)=0.
\end{aligned}
$$
Next using the fact that $\D u={\rm div}{\rm curl}u+\n{\rm div}u$ with $\big(({\rm curl}u)_{ij}=\p_i u_j-\p_j u_i\big)$, we have:
$$
\begin{aligned}
&\rho\p_t v+\rho u\cdot\n v-\mu(\rho){\rm div}{\rm curl}u-2\mu(\rho)\n{\rm div}u-2\n\mu(\rho)\cdot Du-\lambda(\rho)\n {\rm div}u\\
&-\n\lambda(\rho)\,{\rm div}u+\rho^2 f'(\rho)\n {\rm div}u+\rho \n (f'(\rho)\rho)\,{\rm div}u+\rho\n f(\rho)\cdot ^t\n u
+\n P(\rho)=0.
\end{aligned}
$$
It gives in particular:
$$
\begin{aligned}
&\rho\p_t v+\rho u\cdot\n v-\mu(\rho){\rm div}{\rm curl}u-(2\mu(\rho)+\lambda(\rho)-\rho^2 f'(\rho))\n{\rm div}u+(\rho\n f(\rho)-\n\mu(\rho))\cdot ^t\n u\\
&\hspace{5cm}-(\n\lambda(\rho)-\rho\n (f'(\rho)\rho)){\rm div}u-\n\mu(\rho)\cdot\n u+\n P(\rho)=0.
\end{aligned}
$$
In order to introduce a term in $\n\mu(\rho)\cdot{\rm curl}u$ the only possible choice on $f$ corresponds to set $\rho f'(\rho)=2\mu'(\rho)$, we have then:
$$
\begin{aligned}
&\rho\p_t v+\rho u\cdot\n v-\mu(\rho){\rm div}{\rm curl}u-(2\mu(\rho)+\lambda(\rho)-2\rho \mu'(\rho))\n{\rm div}u-\n \mu(\rho)\cdot {\rm curl}u\\
&\hspace{5cm}-(\lambda'(\rho)-2\rho\mu''(\rho)){\rm div}u\,\n \rho+\n P(\rho)=0.
\end{aligned}
$$
In order to "kill" the term in ${\rm div}u$ it is the natural to assume the algebraic relation introduce in \cite{BD}:
$$2\mu(\rho)+\lambda(\rho)=2\rho \mu'(\rho),$$
then:
$$\lambda'(\rho)=2\rho\mu'(\rho),$$
and we deduce that using the fact that ${\rm curl}u={\rm curl}v$ we have the following system with the unknowns $(\rho,v)$:
$$
\begin{cases}
\begin{aligned}
&\p_t\rho-2\D\mu(\rho)+{\rm div}(\rho v)=0,\\
&\rho\p_t v+\rho u\cdot\n v-{\rm div}(\mu(\rho){\rm curl}v)+\n P(\rho)=0.
\end{aligned}
\end{cases}
$$
\hspace{14,6cm}$\blacksquare$
 
\end{document}